\documentclass{amsart}

\usepackage{color}

\newtheorem{propo}{Proposition}[section]
\newtheorem{defi}[propo]{Definition}

\newtheorem{lemma}[propo]{Lemma}
\newtheorem{corol}[propo]{Corollary}

\newtheorem{theo}[propo]{Theorem}

\newcommand{\mar}{\marginpar}

\newcommand{\bl}{\begin{lemma}\label}
\newcommand{\el}{\end{lemma}}

\newcommand{\ld}{,\ldots ,}

\newcommand{\ra}{ \rightarrow }

\newcommand{\se}{ \subseteq }
\newcommand{\lan}{ \langle }
\newcommand{\ran}{ \rangle }

\newcommand{\diag}{\mathop{\rm diag}\nolimits}

\newcommand{\Id}{\mathop{\rm Id}\nolimits}
\newcommand{\Irr}{\mathop{\rm Irr}\nolimits}

\newcommand{\CC}{\mathop{\mathbb C}\nolimits}
\newcommand{\FF}{\mathop{\mathbb F}\nolimits}
\newcommand{\NN}{\mathop{\mathbb N}\nolimits}

\newcommand{\ZZ}{\mathop{\mathbb Z}\nolimits}

\newcommand{\al}{\alpha}

\newcommand{\ep}{\varepsilon}
\newcommand{\lam}{\lambda }

\newcommand{\om}{\omega }

\newcommand{\up}{^{-1}}

\newcommand{\si}{\sigma }

\def\d12{{_{12}}}

\def\acf{{algebraically closed field }}

\def\au{{automorphism }}

\def\ccc{{constituent }}

\def\ei{{eigenvalue }}
\def\eis{{eigenvalues }}

\def\f{{following }}

\newcommand{\med}{\medskip}

\def\hos{{homomorphisms }}

\def\ii{{if and only if }}

\def\ir{{irreducible }}

\def\irt{{irreducible. }}
\def\irr{{irreducible representation }}

\def\st{Suppose that }
\def\hw{highest weight }

\def\itf{{It follows that }}

\def\mult{{multiplicity }}

\def\po{{polynomial }}

\def\rep{{representation }}
\def\reps{{representations }}

\begin{document}

\title[]{Abelian subgroups and semisimple elements in $2$-modular representations of the symplectic group $Sp_{2n}(2)$}

\author[A. Zalesski]{Alexandre Zalesski}

\thanks{e-mail:
 E-mail:   alexandre.zalesski@gmail.com (Alexandre Zalesski)}

\subjclass[2000]{20G05, 20G40} \keywords{Symplectic group, Characteristic $2$, Modular representations, Semisimple elements,
Eigenvalue $1$}
\maketitle

\def\ag{algebraic group }

{\small{\bf Abstract}
We determine the \ir 2-modular \reps of the symplectic group $G=Sp_{2n}(2)$ whose restrictions to every abelian subgroup has a trivial constituent.  A similar  result  is obtained for maximal tori of $G$. There is significant information on
the existence of \ei 1 of elements of $G$ in a given \ir \reps of $G$.
}


\section{Introduction}

The primary motivation of this work lies in the study of \eis of semisimple elements of groups of Lie type
in their \reps over fields of the natural characteristic. This problem is not treatable in full generality,
however, a number of special cases of it  were discussed in literature. In \cite{Z91} there  appeared a number of useful observations, and certain comments are available in surveys \cite{Z88} and \cite{z09}. Occurrance of \ei 1 deserves a particular attention, as this carries some geometrical meaning, and is significant in applications. See \cite{GT}, where one can learn the level of difficulty of the eigenvalue 1 problem. In general, relatively little is known.

Other motivations come from the study of finite groups with disconected prime graphs
 \cite{K,KK}, from  a characterization of groups by the  set of their element orders \cite[Lemma 10]{Zav},  from the study of decomposition numbers \cite{Z16}, from other applications described in \cite{GT} and \cite{z09},
in each case with a specific aspect of the general problem. The complexity of the problem in general justifies considerations of special cases,
one of these is dealt with in the current paper, where we specify the problem to the symplectic group $Sp_{2n}(2)$, $n>1$. Note that the group $GL_n(2)$ is dealt with in \cite{Z17}.



To state our results,
 we need a parametrization of the \ir \reps of $G$. This is available in terms of algebraic groups.
 So we first  recall some fact of algebraic group representation theory, specified here to the group
 $\mathbf{G}=Sp_{2n}(F)$, where $F$ is \acf of characteristic 2.

Note that every 2-modular \irr of $G=Sp_{2n}(2)$ lifts to a \rep of the algebraic group $\mathbf{G}=Sp_{2n}(F)$. The \ir \reps $\phi$ of $\mathbf{G}$ are parameterized by
so called highest weights, equivalently, by the strings of non-negative integers $(a_1\ld a_n)$.
An \irr of $\mathbf{G}$ (and the highest weight of it) is called 2-{\it restricted} if $0\leq  a_1\ld a_n\leq 1$. If $\phi$ is 2-restricted
then the restriction $\phi|_G$ is \ir and all \ir \reps of  $G$ are obtainable in this way. Moreover,  the \ir \reps $\rho$ of $G$ are parameterized by the strings  $\om=(a_1\ld a_n)$ with $0\leq a_1\ld a_n\leq 1$. So we write $\rho_\om$ to specify the parameter.
Set $\om_i=(0\ld 0,1,0\ld 0)$ with 1 at the $i$-th position.

For a group $X$ let $1_X$ be the trivial one-dimensional \rep of $X$.

\begin{theo}\label{ff3}
   Let $G=Sp_{2n}(2)$, and let $\rho$ be an \ir F-representation of G with \hw $ \om=(a_1\ld a_n)$.
Then $1_{A}$ occurs as a constituent of the restriction $\rho|_{A}$ for every abelian subgroup A of G \ii one of the following holds:

$(1)$ $a_n=0$ and either $\sum a_i$ is even or $\sum a_i>2;$

$(2)$ $a_n=1$ and $\sum a_ii\geq 2n$.
\end{theo}

Theorem  \ref{ff3} yields a sufficient condition for an arbitrary element $g\in G$ to have \ei 1
in a given \irr $\rho$. Obtaining a sharp necessary condition for an arbitrary $g$ does not seem to be
a realistic task; in Theorem \ref{si1} we provide a criterion for $\rho(g)$ to have \ei 1 for every $g\in G$.

\begin{defi}\label{d11} Let $n$ be a natural number. Then the Singer height ${\rm Si}(n)$ of n is the maximum number l such that there are natural numbers $n_1\ld n_l $ such that $n_1+\cdots +n_l\leq n$ and the numbers $2^{n_1}+1\ld 2^{n_l}+1$ are coprime to each other.
\end{defi}

Obviously, $1\leq {\rm Si}(m)\leq {\rm Si}(n)$ for $m<n$, but in general the behavior of ${\rm Si}(n) $ is quite irregular.

\begin{theo}\label{si1}
Let $\rho_\om$ be an \irr of $G=Sp_{2n}(2)$ with \hw $\om=(a_1\ld a_n)$.
Then $1$ is an \ei of $\rho(g)$ for every  $g\in G$ \ii one of the following holds:

 $(1)$  $a_n=0$ and $\om\neq \om_i$   for $i$ odd;

 $(2)$   $a_n=1$ and $\sum a_ii\geq n+{\rm Si}(n)$.\end{theo}

Using terminology of  \cite[Definition 1.2]{GT}, a \rep $\rho$ of a group  $G$ is called {\it unisingular} if  $\rho(g)$ has \ei 1 for every   $g\in G$.  So Theorem \ref{si1} classifies \ir   unisingular \reps of $G=Sp_{2n}(2)$ in characteristic 2. Theorem \ref{th1} refines this result for elements of prime-power order.

\begin{theo}\label{th1}
 Let $G=Sp_{2n}(2)$ and let $\rho$ be an \ir $F$-\rep of $ G$ with 
highest weight $\om$. Suppose that $\om\neq \om_i$   for $i$ odd or $i=n$.
If  $g\in G$ is of prime power order  then $1$ is an \ei of $\rho(g)$.
\end{theo}

We have rather precise information on \eis of $\rho_{\om_n}(g)$ for an  arbitrary semisimple element $g\in G$, see Lemma \ref{ff2}, and also of $\rho_{\om}(g)$ for $a_n=1$ (Theorem \ref{fr1}). However,
the occurrence of \ei 1 of $\rho_{\om_i}(g)$ with $i<n$ odd for an arbitrary $g\in G$ does not reveal
good regularities, and we have no uniform result even for $p$-elements $g$.
If $a_n=0$ then  $\rho(g)$ has \ei 1 for every $g\in G$ whenever so has $\rho(t)$ for $t\in G$ of order $2^n+1$. Similar "testing" elements exist also when  $a_n=1$, they are constructed in terms of ${\rm Si}(n)$,  see Proposition \ref{p99}. 

Another  test is Proposition \ref{p88}:  for $g\in G$ given, 1 is an \ei of $\rho(g)$ if 1 is an \ei of  $\rho_{\om_1}(g)$ and $\rho_{\om_1}(g)$.


As a by-product we obtain the \f result on the occurrence of weight 0 in an \irr  $\rho_{\om}$ of the algebraic group $Sp_{2n}(\overline{\FF}_2)$.

 \begin{theo}\label{m22a}
   Let $\om=(a_1\ld a_n)$ be a $2$-restricted dominant weight of ${\mathbf G}=Sp_{2n}(\overline{\FF}_2)$, and let $\om'=(a_1\ld a_{n-1},0)$. Suppose that $a_n=1$. Then the following statements are equivalent:

$(1)$ V has weight $0;$

$(2)$ $\sum a_i i$ is even and greater than $2n-1;$

$(3)$ $\om_n\prec\om'.$
\end{theo}

Note that, in view of \cite[Theorem 15]{z09}, the case considered in Theorem \ref{m22a}, was the only one where the occurrence of weight 0 in an infinitesimally \irr of a simple \ag was not known.

Another result of independent interest is

\begin{theo}\label{th7a}
 Let F be an \acf of any characteristic, $\mathbf{ G}=Sp_{2n}(F)$ and
$n=n_1+\cdots+n_k$, where $n_1\ld n_k$ are positive integer. Let $\mathbf{ H}=\mathbf{ H}_1\times \cdots\times \mathbf{ H}_k$, where $\mathbf{ H}_i\cong Sp_{2n_i}(F)$ for $i=1\ld k$.
Let $\rho$ be an  \irr of G with p-restricted \hw $\om=(a_1\ld a_n)$. Then the \f are equivalent:

$(1)$ $\tau(\mathbf{ H}_i)=1$ for some $i\in\{1\ld k\}$ whenever $\tau$ is a composition factor   of $\rho|_\mathbf{ H};$

$(2)$ $\sum a_ii<k$.\end{theo}

Note that the case with $k=2$ in Theorem \ref{th7a} follows from the main result of \cite{SZ}.

 Recall that every element of $G$ of odd order lies in a maximal torus of $G$. So one can ask when the trivial \rep $1_T$ of a maximal torus $T$ of $G$ occurs as a constituent of the restriction of $\rho_\om$ to $T$. The answer is similar to the result of Theorem \ref{ff3}, see Theorem \ref{ee3}.

In \cite{Z17} we classify unisingular \ir \reps of $H=GL_n(2)=SL_n(2)$.   For some applications one needs a more precise information for specific elements $h\in H$. We apply the above results to deal with the real elements of  $H$. Recall that an element $x$ of a group $X$ is called {\it real} if $x$ is conjugate to $x\up$.

\begin{theo}\label{th6}
Let $H=SL_{2m}(2)$ or $SU_{2m}(2)$, and $\rho$  an \ir $2$-modular \rep of $H.$ Let $h\in H$ be real. Then $\rho(h)$ has \ei $1$ unless, possibly,
$\rho$ is an odd exterior power of the natural (that is, defining) \rep of $H$.\end{theo}

 Note that $h\in SL_{2m}(q)$ is real whenever $|h|$ divides $q^i+1$ for some $i$  (Lemma \ref{yz2}),
 and $h\in SU_{2m}(q)$ is real whenever $|h|$ divides $q^i-1$ for some odd $i$  (Lemma \ref{ru1}).

Observe that if $H=SL_n(F)$ with $n$ odd and $G=Sp_{n-1}(F)$ for any \acf $F$ then
$\rho|_{G}$ has weight 0 for an arbitrary \irr $\rho$ of $H$ (Lemma \ref{hg5}(3)); therefore $\rho(g)$
has \ei 1 whenever $g\in H$ is real (Lemma \ref{hg5}).

 Kondrat'ev \cite{K}    raises the problem of
determining the \ir modular \reps of simple groups $G$ over a prime field such that 1 is not an \ei
of $\rho(g)$ for $g$ to be a prime order. (In \cite{KOS}
the problem was specified to \ir elements  $g\in G=SL_n(q)$ such that $|g|$ is a prime dividing $(q^n-1)/(q-1)$ and $q$ a $p$-power.) So Theorems \ref{th1} and \ref{th6} contribute to this problem.

Other results of  the paper refine the above results for some special cases. We mention the following:
 
\begin{theo}\label{th2}\mar{th2}  Let $\phi$ be an \ir $F$-\rep of $ \mathbf{G}$ with highest weight $\om=(a_1\ld a_n)$.
Let s be a Singer cycle of $G$. Then  $1$ is not an \ei of $\phi(s)$ \ii 
 $\om=2^k\om_n$ or  $\om=2^k\om_i$ for some odd $i\in\{1\ld n-1\}$ and an integer $k\geq 0$.
\end{theo}

An element of order $2^n+1$ in $G=Sp_{2n}(2)$ is called a {\it Singer cycle} \cite{Hu}. This generates a maximal torus of $G$ and plays a role in some applications.

\med
{\it Notation}\,\, $\CC$ is the field of complex numbers, $\ZZ$ is the ring of integers and $\NN$ the set of natural numbers. For  integers $a,b>0$ we write $(a,b)$ for the greatest common divisor of $a,b$, and $a|b$ means that $b$ is a multiple of $a$. For sets $A,B$ in a $\ZZ$-lattice we write $A+B$ for the set $\{a+b:a\in A,b\in B\}$. A diagonal $(n\times n)$-matrix with subsequent entries $x_1\ld x_n$ is denoted by $\diag(x_1\ld x_n)$.  By $\FF_q$ we denote the finite field of $q$ elements and by $\overline{\FF}_q$ the algebraic closure of $\FF_q$.
All \reps below are over an \acf $F$ of characteristic 2 (unless otherwise is stated explicitly). $F^\times$ is the multiplicative group of $F$.
If $\phi$ is a \rep of a group $X$ and $Y\subset X$ is a subgroup then $\phi|_Y$ means the restriction of $\phi$ to $Y$. The tensor product of \reps $\phi,\psi$ is denoted by $\phi\otimes\psi$.

Notation and terminology for weight systems are introduced in Section 2.1, for those for algebraic groups
 see Sections 2.2 and 2.3.

    \section{Preliminaries}

\subsection{The weight system of type $C_n$}

For a precise definition of weight system of type $C_n$ see \cite{Bo}. We recall
here some notions of the  theory.

A {\it weight lattice} $\Omega$ is formed by
 strings of integers $(a_1\ld a_{n})$, which are   called {\it weights of} $\Omega$, and those with non-negative
entries $a_1\ld a_{n}$ are called {\it dominant weights}.
The subset of dominant weights is denoted by $\Omega^+$.
The weights  $\om_i=(0\ld 0,1,0\ld 0)$ $(1\leq i\leq n)$, where $1$ occupies the $i$-th position, are
called {\it fundamental}, and the weight $(0\ld 0)$ is called the zero weight. We usually simplify this by denoting this by 0. So arbitrary weights can be written as $\sum_{i=1}^{n} a_i\om_i$ with $a_1\ld a_n\in \ZZ$.  For an integer $m>1$ a dominant weight
$(a_1\ld a_{n})$ or $\sum_{i=1}^{n} a_i\om_i$ is called $m$-{\it restricted} if $a_i<m$ for $i=1\ld n$. (We use this mainly for $m=2.$) We denote by $\Omega'$ the subset of weights $(a_1\ld a_{n-1},0)$.

The theory of weight system of type $C_n$ singles out a sublattice $\mathcal{R}$   of index $2$ in $\Omega$ called the {\it root lattice}. In addition, one singles out some elements   $\al_1\ld \al_{n}\in \mathcal{R}$ called {\it simple roots}, which form a $\ZZ$-basis of $\mathcal{R}$. Denote by $\mathcal{R}^+$ the set of non-negative
linear combinations of   simple roots.  For  weights $\om,\mu$ one writes $\mu \preceq \om$
if $\om-\mu\in \mathcal{R}^+$; if additionally $\om\neq \mu$ we write  $\mu \prec \om$. If $\mu \prec \om $ are dominant then we say that $\mu$ is a subdominant weight for $\om$.
The   weights $\om\in \mathcal{R}$ are called {\it radical}. We write $0\prec \mu$ to state that the weight 0 is a subdomininant weight of $\mu$.

The  theory of weights can be viewed as a part of linear algebra, and it is somehow independent from the theory of algebraic groups, see for instance \cite{Bo}, where detailed data on relations between
 fundamental weights and simple roots are tabulated in \cite[Planche III]{Bo}.

In our analysis an essential role is plaid by another basis $\ep_1\ld \ep_n$ of $\Omega$; specifically $\ep_1=\om_1$ and
$\ep_i=\om_i-\om_{i-1}$ for $i=2\ld n$.

We have $\om_i=\ep_1+\cdots +\ep_i$ for $i=1\ld n$, and $\om_i$ is radical \ii $i$ is even. The simple roots are $\al_i=\ep_i-\ep_{i+1}$
for $i<n$, and $\al_n=2\ep_n$. Note that $\al_i=2\om_i-\om_{i-1}-\om_{i+1}$ for $i<n$, and $\al_1=2\om_1-\om_{2}$, $\al_n=2\om_n-2\om_{n-1}$. 

The Weyl group $W$ of $\mathbf{G}$ (or the weight system $\Omega$) acts on  $\Omega$ as a group of linear transformations.
Specifically, $W$ transitively permutes $\pm\ep_1\ld \pm\ep_n$. In fact, $W$ has an elementary abelian normal subgroup $W_1$
which acts on $\pm\ep_1\ld \pm\ep_n$ by changing the signs, and $W$ is a semidirect product $ W_1\cdot S_n$, where $S_n$ is the symmetric group; the latter acts on $\ep_1\ld \ep_n$ by permutations.
It follows that the $W$-orbit of $\om_i=\ep_1+\cdots+\ep_i$ consists of weights $\pm \ep_{j_1}\pm\cdots\pm \ep_{j_i}$,
where $\ep_{j_1}\ld  \ep_{j_i}\in\{\ep_1\ld \ep_n\}$ and $j_1<\cdots<j_i$.

For a weight $\mu=\sum a_i\om_i$ set  $\delta(\mu)=\sum a_ii$ and $\gamma(\mu)=\sum a_i$.
 
\bl{sd3}
 $(1)$ $\delta(\al_i)=0$ for $i=1\ld n-1$, and $\delta(\al_n)=2$.

$(2)$ $\gamma(\al_i)=0$ for $i=2\ld n-1$, and $\gamma(\al_1)=1=\gamma(\al_n)$.

$(3)$ Let $\om,\nu\in\Omega$. If $\om>\mu$ then $\delta(\om)\geq \delta(\mu)$ and  $\delta(\om)- \delta(\mu)$
is even. In particular, $\om\in {\mathcal R}$  \ii $\delta(\om)$ is even.
\el

Proof. (1) and (2) follow straightforwardly from the above expressions of the simple roots in terms of the dominant weights,
and (3) follows from (1) and the definition of the partial order $\om>\mu$.

\bl{sd4}
Let $\om$ be a dominant weight of $C_n$.

$(1)$ If\,  $0 \prec\om$  then $\om_2\preceq \om.$

$(2)$ If\,  $\om$ is not radical then $\om_1\preceq\om.$

$(3)$
 Let 
 $\om=\sum a_i\om_i$ 
and $c=\sum a_i$. Then $ c\om_1\preceq \om $ or $(c-1)\om_1\prec (c-1)\om_1+\om_2\prec \om.$
In addition, if $\om_1\prec\om$ and $c>2$ then $3\om_1 \preceq\om$.

$(4)$ $\om_{i-2}\prec \om_i$  for $i=3\ld n$. Moreover,  $\om\prec\om_i$ \ii $i<j$ and $j-i$ is odd.

$(5)$ Suppose that $\om\neq \om_i$ for $1\leq i\leq n$. If $\om_1\prec\om$ then $\om_1+\om_2\preceq \om$; if $0\prec\om$ then $2\om_1\preceq\om.$

$(6)$ 
Let\,  $l=\delta(\om)$, and $l=j+kn$ with $0\leq j<n$. If $j=0$ then $k\om_n\preceq \om$, otherwise $\om_i+k\om_n\preceq \om$.
In addition, if $k>0$ then $\om_n\preceq \om$ if $l-n$ is even, and $\om_1+\om_n\preceq \om$ if\, $l-n$ is odd.
 \el

Proof. The statements (1),(2) are well known.

(3) Set $c_1=\sum_{i~odd} a_i$,  $c_2=\sum_{i~even} a_i$.
By (1),(2) we have $\om_1\preceq\om_i$ for $i$ odd and $\om_2\preceq\om_i$ for $i$ even. Therefore,
$ c_1\om_1+c_2\om_2\preceq \om$. As $ 2\om_1\prec 2\om_2$,  the result follows if $c_2$ is even, otherwise
we have $(c-1)\om_1\prec (c-1)\om_1+\om_2\prec $. If  $\om_1\prec\om$ then $c$ is odd in the former case, whence the claim,
and $c$ is even in the latter case. Then $c\geq 4$, and the result follows.

(4) The first claim is obvious as $\om_i-\om_{i-2}=\ep_i+\ep_{i-1}$ is a positive root \cite[Planche III]{Bo}. For the second one
follows from Lemma \ref{sd3}(3).

(5) By (3), $ c\om_1\preceq \om$ or $ (c-1)\om_1+\om_2\preceq\om$. If $\om_1\prec\om$ then $c$ is odd in the former case
and even in the latter one. As $(c-2)\om_1\prec c\om_1$ for $c>2$, we have $ \om_1+\om_2\prec 3\om_1\preceq \om$ in the former case, and $\om_1+\om_2\preceq \om$ in the latter case. If $0\prec\om$ then $c$ is even in both the cases, whence
$ 2\om_1\preceq\om$ or $2\om_1\prec 2\om_1+\om_2\preceq \om$, as required.

(6) It suffices to prove this for $\om $ minimal, in the sense  that there is no dominant weight  $\mu\prec\om$ such that $\delta(\mu)= \delta(\om)$. As $\delta(\al_i)=0$ for $i<n$, it follows that $\om-\al_i$ is not a dominant weight for every $i=1\ld n-1$. This implies $a_i\leq 1$ for $i=1\ld n-1$.

Suppose that there are $i<j<n$ such that $a_i=a_j=1$. As $\om_i+\om_j=\om_{i-1}+\om_{j+1}+\al_i+\cdots +\al_j$, we have $\om_{i-1}+\om_{j+1}\prec\om_{i}+\om_{j}$ and $\delta(\om_{i}+\om_{j}) =\delta(\om_{i-1}+\om_{j+1})$ (because $j<n$). This is a contradiction. So there is at most one $i<n$ such that $a_i>0$, so $\om=k\om_n$ or $\om_i+k\om_n$, as required.

This also implies the additional statement if $l\leq n+1$. Suppose $l>n+1$. Now choose $\om$ to be minimal in the sense  that there is no dominant weight  $\mu\prec\om$ such that $\delta(\mu)>n+1$. If
 $a_n> 1$ then $\delta(\om-\al_n)=l-2$,  contradicting the minimality of $\om$  unless
$l-2\leq n+1$, that is, $l\leq n+3$. In this case $n+3\geq l=\delta(\om)\geq \delta(2\om_n)=2n$, whence $n\leq3$.
If $n=3$ then $\om=2\om_3$ and $\om-\al_3=2\om_2>\om_1+\om_3$ as required. If $n=2$ then $\om=\om_1+2\om_2$ or $2\om_2$, and
the lemma is true in these cases as  $\om_2\prec 2\om_1\prec2\om_2$.

Let $a_n= 1$. Then $\om=\om_i+\om_n$ and $0\prec\om_i$ if $i$ is even, otherwise $\om_1\prec\om$, whence the result. $\Box$

\subsection{Abelian subgroups and maximal tori: generalities}

Let $V$ be the natural $\FF_2G$-module and $A$ an abelian  subgroup of $G$ of odd order. By Maschke's theorem,
$V=U_1\oplus\cdots \oplus U_m$, where $U_1\ld U_m$ are irreducible $F_2A$-modules. Each $U_i$, $i=1\ld m$,  is either non-degenerate or totally isotropic; in the latter case there is a unique $j$ such that $U_i+U_j$ is non-degenerate.
Moreover, $U_j$ is dual to $U_i$ as $F_2A$-modules. We can rearrange the above decomposition and write
$V=V_1\oplus \cdots\oplus V_k$, where each $V_i$ is non-degenerate and either \ir or a sum of two \ir submodules dual to each other.
Let $2n_i=\dim V_i$ and let $A_i$ be the restriction of $A$ to $V_i$. Then $A_i\subset Sp(V_i)$ is a cyclic group and  $|A_i|$ divides $2^{n_i}+1$ if $V_i$ is irreducible, otherwise
$|A_i|$ divides $2^{n_i}-1$.  In addition, $n_1+\cdots +n_k=n$. If $A$ is a maximal abelian subgroup of odd order then $A\cong A_1\times \cdots\times A_k$
and $|A_i|=2^{n_i}+1$ or $2^{n_i}-1$, in the latter case $n_i>2$. If $k=1$ and $|A_1|=2^{n}+1$ then $A$ is called  a {\it Singer subgroup of} $G$, and the generators of $A$ are called {\it Singer cycles}.

In the theory of finite reductive groups an important role belongs to maximal tori, which are abelian groups
of order coprime to the characteristic of the ground field. Every maximal abelian subgroup of odd order is a maximal torus, but not conversely. The definition and general properties of maximal tori is available in
\cite{DM} and \cite{Ca}, for specific properties of maximal tori in classical groups  used in this paper one can
consult with \cite{Z16}, \cite{Z14} or \cite{HZ09}. The maximal tori in $G=Sp_{2n}(2)$ are in bijective correspondence with the decompositions  $V=V_1\oplus \cdots\oplus V_k$ defined above, and the conjugacy classes of maximal tori in $G$ are in bijective correspondence with the $G$-orbits of the decompositions. These can be
parameterized by the strings $(\pm n_1\ld \pm n_k)$ up to the ordering, the sign minus is chosen whenever the torus acts
on $V_i$ irreducibly ($i=1\ld k$). The string $(1\ld 1)$ corresponds to the torus of one element $\Id$,
the string $(-1\ld -1)$ labels the torus denoted below by $T^\#$ (so $|T^\#|=3^n$). A maximal torus of order $2^n+1$ is called a {\it Singer torus}; it corresponds to $(-n)$. A maximal torus corresponding to $(n)$ is a cyclic group of order $2^n-1$.
Every maximal torus $T$ of $G$ can be written as $T=T_1\times\cdots \times T_k$, where $|T_i|=2^{n_i}+1$
if the sign of $n_i$ is minus, otherwise $|T_i|=2^{n_i}-1$.  (The above condition $n_i>2$ is irrelevant now.)
The decomposition is unique up to the ordering of $T_1, \ldots, T_k$, and we call it the {\it cyclic decomposition}
(reflecting the fact that $T_1,\ldots, T_k$ are cyclic groups). Observe that the multiples $T_i$ with $|T_i|=1$ ($1\leq i\leq k$) are not dropped from the expression $T=T_1\times\cdots \times T_k$.

Sometimes it is convenient to use
a pair of strings $(n_1\ld n_k)$, $(\eta_1\ld \eta_k)$, where $\eta_1\ld \eta_k\in\{1,-1\}$. Then
$|T_i|=2^{n_i}-\eta_i$ for $1\leq i\leq k$. \itf every maximal torus is contained in a direct product
$G_T=Sp_{2n_1}(2)\times \cdots \times Sp_{2n_k}(2), $ where $T_i$ is a maximal torus of $Sp_{2n_i}(2)$ for
$i=1\ld n$.
In addition, $T_i$ is a Singer torus in $Sp_{2n_i}(2)$ \ii  $\eta_i=-1$,
otherwise $T_i$ is of order $2^{n_i}-1$ and reducible in $Sp_{2n_i}(2)$. The number of indices  $i$
such that $\eta_i=-1$ is here called the {\it Singer index of} $T$.

Let $T=\lan t\ran$ 
be a maximal  torus of order $2^n+1$ and $2^n-1$,
respectively. (These are known to be unique up to conjugation, see for instance \cite[Lemma 2.1]{Z14} and \cite[Lemma 7.1]{EZ}.) Then $t$ is conjugate in $\mathbf{G}$ to an element $d\in {\mathbf T}$ such that $\ep_i(d)=\zeta^{2^{i-1}}$ for $i=1\ld n$, where $\zeta$ is a primitive $|t|$-root of unity. (This is well known and explained in detail in \cite{Z16}.) Under a certain basis of $V,$ the matrix of $J$ on $V$ is diagonal of shape
\begin{equation}d:=\diag(J,J\up )\in \mathbf{G}=Sp_{2n}(\overline{\FF}_{2}),~{\rm where}~ J=\diag(\zeta,\zeta^2\ld \zeta^{2^{n-1}})
.\end{equation}
Therefore, $\ep_i(d)=\zeta^{2^{i-1}}$ for $i=1\ld n$. This implies

\begin{equation}\label{eq2} 2\ep_i(d)=\ep_{i+1}(d)~ {\rm for} ~i=1\ld n-1~{\rm and }~ 2\ep_n(d)=\zeta^{2^n}=\begin{cases}-\ep_{1}(d)&if ~|t|=2^n+1\\ ~\ep_{1}(d)&if ~|t|=2^n-1.\end{cases} \end{equation}

This can be extended to arbitrary maximal tori. Let $V=\oplus V_i$ and let $T_i$ be a maximal torus in $Sp(V_i)\cong Sp_{2n_i}(2)$,
so $T$ is conjugate to the group  $D:=\{\diag(J^{m_i}_1\ld J^{m_k}_k,J^{-m_k}_k\ld J^{-m_i}_1\}$, where $J_i$ $(i=1\ld k)$ is constructed as above,  $1\leq m_i\leq |T_i|$.
We call $D$ a {\it canonical form of} $T$. Note that there is a conceptual way to construct $D$ for arbitrary
finite group of Lie type, see the proof of \cite[Proposition 3.3.6]{Ca}.

\bl{ss1}
 \cite[p. 209, Corollary 5.10(a)]{SS} Let $G=Sp_{2n}(2)$ and $A\subset G$ an abelian subgroup of odd order. Then A is contained in a maximal torus of G. In particular, if A is a maximal
abelian  subgroup of odd order then A is a maximal torus of G.
\el

The \f lemma is trivial, but we state it explicitly for reader's convenience.

\bl{aa1}\mar{aa1} Let A be an abelian subgroup of G, and $\rho$ a \rep of G over a field of characteristic $2$.
Let B be a maximal subgroup of A of odd order.
Then $1_A$ is a \ccc of $\rho|_A$ \ii $1_B$ is a \ccc of $\rho|_B$.
\el

\subsection{Representations and their weights}

The above parametrization of the \ir \reps of $G=Sp_{2n}(2)$
is based on Steinberg's theorem \cite[\S 2.7]{Hum} saying that the \ir
\reps of $Sp_{2n}(2)$ are in bijection with the  $2$-restricted dominant weights of $\mathbf{G}$.
In fact, the \rep theory of $\mathbf{G}$ plays a significant role in proving the above result. The \rep theory of the algebraic group
$Sp_{2n}(F)$ is described on the language of the theory of the weight system of type $C_n$.

Let $F$ be an \acf of characteristic $p>0$. The group $\mathbf{G}=Sp_{2n}(F)$ is defined as the group of $(2n\times 2n)$-matrices over $F$ preserving a non-degenerate alternating form $(*,*)$ on $V=F^{2n}$. A basis $e_1\ld e_{2n}$ of $V$ is called a Witt basis if $(e_i,e_{2n+1-i})=1$ and
$(e_i,e_j)=0$ if $j\neq 2n+1-i$ for $i=1\ld n$ and $1\leq j\leq 2n$.  The subgroup  $\mathbf{T}$ of diagonal matrices in  $\mathbf{G}$
under a Witt basis is of shape $\diag(t_1\ld t_n,t_n^{-1}\ld t_1\up)$ ($t_1\ld t_n\in F^\times $).
We refer to this group as a reference torus, and define a maximal tori in  $\mathbf{G}$ as those conjugate to  $\mathbf{T}$
in $ \mathbf{G}$. The weights are defined as rational \hos $\mathbf{T}\ra F^\times$, and $\ep_i$ are defined as those sending
 $\diag(t_1\ld t_n,t_n^{-1}\ld t_1\up)$  to $t_i$.

The Weyl group $W$ of $\mathbf{G}$ is defined as $W:=N_{\mathbf{G}}(\mathbf{T})/\mathbf{T}$. The conjugation action of $N_{\mathbf{G}}(\mathbf{T})$ on $\mathbf{T}$ yields the action of $W$ on $\Omega={\rm Hom}(\mathbf{T},F^\times)$.  In this action $W$ preserves the set $\{\pm\ep_1\ld \pm\ep_n\}$ and acts
on this set transitively.
Note that the $W$-orbit of any weight $\om\in\Omega$ contains a unique dominant weight.

For a $\mathbf{G}$-module $V$ or a \rep $\rho$ of $\mathbf{G}$ we  denote by $\Omega(V)$, $\Omega(\rho)$
the set of weights of $V, \rho$, respectively. Then $\Omega(V)$ is $W$-invariant, in particular, $\Omega(V)=-\Omega(V)$.

Let 
$\mu$ be a dominant weight and $W$ the Weyl group of $\mathbf{G}$. 
Denote by $\chi_\mu$ the sum  of distinct  weights $w(\mu)$ with $w\in W$.
For a semisimple element $g\in G$ let $g'$ be a conjugate of $g$ in $ \mathbf{T}$ (it is well known that every semisimple element of $\mathbf{G}$ is conjugate to one in $\mathbf{T}$). Then $g'$ is not unique but the value $\chi_\mu(g')$ does not depend from the choice of $g'\in\mathbf{T};$ so $\chi_\mu$ is a well defined function
  on the semisimple elements of $G$ called an {\it orbit character}. (This is well known to be a generalized Brauer character of $G$, see \cite[\S 5.7]{Hum}.)

If $F=\overline{\FF}_p$ and $\beta_\om$ is the Brauer character of  $V_\om$ then
  $\beta_\om=\sum m_\mu\chi_\mu$, where $\mu$ runs over the dominant weights of  $V_\om$ and $m_\mu$ is the \mult of $\mu$ in $V_\om$, see \cite[\S 2.2 and \S 5.6]{Hum} for details. Therefore, $\beta_\om|_G=\sum m_\mu\chi_\mu|_G$; this is used below for computations of \ei 1 of elements $g\in G$ on $V_\om$.

 \bl{pr4}
   Let F be an \acf of characteristic $p>0$ and $\mathbf{G}=SL_{n+1}(F)$ or $Sp_{2n}(F)$.
 Let V be an \ir $\mathbf{G}$-module with $p$-restricted \hw  $\om=(a_1\ld a_n)$.

$(1)$  The weights of V are the same as those of the \ir L-module $M$ of the Lie algebra over the complex numbers of the same type as $\mathbf{G}$ with
with \hw $(a_1\ld a_n)$   (but the multiplicities of the same weight in V and M may differ), unless
possibly $p=2$, $\mathbf{G}= Sp_{2n}(F)$ and $a_n=1$.  

$(2)$ Let $\mu\prec \om $ be dominant weights. Then $\mu\in\Omega(V)$. Moreover, all weights of $\Omega(V_\mu)\subset\Omega(V_\om)$ provided $a_n=0$ for $p=2$ and $\mathbf{G}=Sp_{2n}(F)$.

$(3)$ Let $\lam$ be a dominant weights for $\mathbf{G}$ such that $\lam+\om$ is p-restricted.
Then $\Omega(V_{\lam+\om})=\Omega(V_{\lam})+\Omega(V_{\om})$.
\el

Proof.  For (1) see \cite[Theorem 15]{z09}. (2) follows from (1) and  \cite[Ch. VIII, Prop 5]{Bo8}.
(3) Clearly, either $\lam\in\Omega'$ or $\om\in\Omega'$ (or both). By swapping $\lam,\om$ we can assume that $\lam\in\Omega$. (i) Suppose first that we are not in the exceptional situation of (1). Then the result follows from (1) and  \cite[Ch. VIII, Prop 10]{Bo8}. (ii) Let
$p=2$, $\mathbf{G}= Sp_{2n}(F)$ and $a_n=1$. Set $\mu'=\om-\om_n$.
Then $V_{\lam+\om}=V_{\lam+\mu'+\om_n}=V_{\lam+\mu'}\otimes V_{\om_n}$, and then
$\Omega(V_{\lam+\om})=\Omega(V_{\lam+\mu'}\otimes V_{\om_n})=\Omega(V_{\lam+\mu'})+\Omega( V_{\om_n})$. By (i), $\Omega(V_{\lam+\mu'})=\Omega(V_{\lam})+\Omega(V_{\mu'})$, so $\Omega(V_{\lam+\om})=\Omega(V_{\lam})+\Omega(V_{\mu'})+\Omega( V_{\om_n})$. As $V_{\mu'}\times V_{\om_n}=V_{\mu'+\om_n}$, we have $\Omega(V_{\mu'})+\Omega( V_{\om_n})=\Omega(V_{\mu'+\om_n})$,
whence the result. $\Box$

\med
If $\om\notin\Omega'$ then the  statement (2) of Lemma \ref{pr4} is not valid anymore;  Lemma \ref{t29} and Theorem \ref{m22} below describe some special cases.

\bl{t29}
  Let $ \om=(a_1\ld a_n)$ be a $2$-restricted  dominant weight of $\mathbf{G}=Sp_{2n}(F)$. 
Suppose that  $a_n=1$ and $\delta(\om)=\sum  a_ii\geq 2n$.

$(1)$\, If\, $\om$ is radical then $0$ is a weight of  $V_\om;$

$(2)$\, if\, $\om$ is not radical then $\om_1,\om_1+\om_2,3\om_1$ are weights of $V_\om$ as well as $\om_3$ for $n>2$.
In addition, all weights of $V_{\om_1+\om_2}$ as weights of $V_{\om}$.\el

Proof. Let $\om'=\om-\om_n$. Then $V_\om=V_{\om'}\otimes V_{\om_n}$. Then $\om'\in\Omega'$, and hence $\mu\prec\om'$ for $\mu$ dominant implies  $\mu$ to be a weight of  $V_{\om'}$ (Lemma \ref{pr4}(2)).

(1) As $\delta(\om')\geq n$ and $\om$ is radical, we have $\om'-\om_n\in\mathcal{R}$, and then $\om_n\preceq \om'$ by Lemma \ref{sd4}(6).
So $\om_n$ is a weight of  $V_{\om'}$. As $\om_n$ is a weight of  $V_{\om_n}$, the claim follows.

(2)  As $\om$ is not radical, we have $\om'-\om_n\notin\mathcal{R}$, so $\om_n$ is not a weight of $V_{\om'}$.
 By Lemma \ref{sd4}(6),  $\om_1+\om_n\preceq \om'$ (as  $\delta(\om')\geq n+1$). So
 $\om_1+\om_n\in \Omega(V_{\om'})$  by Lemma \ref{pr4}, and hence $\om_1\in \Omega(V_{\om})$.

We have $\om_1+\om_n=2\ep_1+\ep_2+\cdots +\ep_n$ and $\om_n=\ep_1+\ep_2+\cdots +\ep_n$. As $\mu=\ep_1-\ep_2-\cdots -\ep_n$ is a weight of $V_{\om_n}$, we conclude that  $3\om_1=\om_1+\om_n+\mu$
is a weight of  $V_{\om}$.

Similarly, $\nu=\ep_1+2\ep_2+\ep_3+\cdots +\ep_n\in\Omega(V_{\om'})$, so
  $\nu+\mu=2\ep_1+\ep_2=\om_1+\om_2 \in\Omega(V_{\om})$.

 Finally, if $n>2$ then  $\om_1+\om_{n-2}\prec \om_1+\om_n$, and hence $\om_1+\om_{n-2}\in \Omega(V_{\om'})$  by Lemma \ref{pr4}.
As $\om_1+\om_{n-2}=2\ep_1+\ep_2+\cdots +\ep_{n-2}$, it follows that $\ep_1-\ep_{n-1}-\ep_n\in\Omega(V_{\om})$. This weight is in the $W$-orbit of $\om_3=\ep_1+\ep_{2}+\ep_3$, so $\om_3 \in\Omega(V_{\om})$.

The additional statement follows as the weights $\om_1$, and $\om_3$ for $n>2$, are the only subdominant weights of  $V_{\om_1+\om_2}$.
$\Box$

  \begin{theo}\label{m22}
   Let $\om=(a_1\ld a_n)$ be a $2$-restricted weight and $a_n=1$. Then
  $V=V_{\om'}\otimes V_{\om_n}$, where $\om'=\om-\om_n$, and 
 the following statements are equivalent:

$(1)$\, V has weight $0;$

$(2)$\,\, $\om_n\prec\om';$

$(3)$ \,  $\delta(\om)$ is even and greater than $2n-1$. 
\end{theo}

Proof. Note that  $\Omega(V)=\Omega(V_{\om'})+\Omega(V_{\om_n})$.

$(2)\ra (1)$. By Lemma \ref{pr4}(1), $\om_n$ is a weight of $V_{\om'}$, and hence so is $-\om_n$. Whence the claim.

$(1)\ra (2)$.  The weights of $V_{\om_n}$ are $W$-conjugate. \itf if a weight $\mu\in\Omega(V_{\om_n})$ then 
$\pm\om_n$ is a weight of $V_{\om'}$. Then $\om_n$ is a weight of  $V_{\om'}$ too, whence $\om_n\prec\om'$.

$(2)\ra (3)$. We have $\om'-\om_n\in\mathcal{R}$ as $\om_n\prec \om'$, and hence $\om'+\om_n\in\mathcal{R}$ as $2\om_n\in\mathcal{R}$. So the result follows from Lemma \ref{sd3}(3).

$(3)\ra (2)$ As $\delta(\om)$ is even, $\delta(\om')-n=\delta(\om)-2n$ is even too.
 Since $\om'$ is dominant and $\om'\neq \om_n$ by assumption,
the result follows from Lemma  \ref{sd4}(6) applied to $\om'$.  $\Box$

\med
Remark. Under the assumption of Theorem \ref{m22}, if $V$ has weight 0 then $\om_i\in\Omega(V)$ for every $i\leq n$ even.
Indeed, by Lemma \ref{pr4}.
$\om_{n-i}\in\Omega(V_{\om'})$ for $i<n$ even. So $\mu:=-\ep_{i+1}- \cdots -\ep_n\in\Omega(V_{\om'})$
so $\mu+\om_n=\ep_1+\cdots+\ep_i\in\Omega(V_{\om})$.


\section{Symplectic groups over a field of two elements}
\def\hw{highest weight }

In this section $G=Sp_{2n}(2)$ and ${\mathbf G}=Sp_{2n}(F)$, where $n>1$ and $F$ is an \acf of characteristic 2. We fix a reference torus ${\mathbf T}$ of ${\mathbf G}$ which determine the weight system of ${\mathbf G}$.
Recall that the conjugacy classes of  maximal tori of $G$ are in a canonical bijection with
the conjugacy classes of the Weyl group $W$ of ${\mathbf G}$.  So one can write $T_w$ for a maximal torus of $G$ labeled by the conjugacy class of $w\in W$. When we consider properties of $T$ in a \rep of ${\mathbf G}$
we choose a canonical conjugate $D\subset  {\mathbf G}$ as defined in Section 2.2. Let $T$
be a maximal torus of $G$. Then $T=T_1\times\cdots\times T_k$ for some $k\leq n$, and let $(\pm n_1\ld \pm n_k)$ be the label of $T$, see Section 2.2. Abusing notation, we often assume that $D=T$.

The \f lemma is a special case of \cite[Lemma 3.9]{Z16}.

\bl{016}
Let $\chi_{\om_i}$ be the orbit character of $G$ and $T=T_w$ a maximal torus
labeled by an element $w\in W$. Then $1_T$ is a constituent of $\chi_{\om_i}|_T$ \ii w is conjugate in W to an element of $W_i$, where $W_i$ is the stabilizer of $\ep_1+\cdots +\ep_i$ in W. \el

Proof. In \cite[Lemma 3.9]{Z16} the result is stated in a more precise form, specifically, it states that the \mult of $1_T$ in $\chi_{\om_i}|_T$ equals the value of the induced character $1_{W_i}^W$ at $w$, where $W_i$ is the stabilizer of $\om_i=\ep_1+\cdots +\ep_j$ in $W$.  It is well known that
this is non-zero \ii $w$ is conjugate to an element of $W_i$. $\Box$

\bl{op1}
Let $s$ be a Singer cycle of G. Then s has no fixed point on $V_{\om_i}$ \ii $i=n$ or $i<n$ odd.
\el

Proof. Set $S=\lan s \ran$. If $i<n$ is even then $V_{\om_i}$ has weight 0, whence the claim for this case. Let $i<n$ be odd.
By Lemma \ref{sd4}(4),  $\om_j$ with $j<i<n$ odd are  the only subdominant weights of $\om_i$.
In addition, the weights of $V_{\om_n}$ are in the $W$-orbit of $\om_n. $  Therefore, it suffices to show that $1_S$ is not a constituent of the
orbit character $\chi_{\om_j}$ for $j$ odd and $j=n$. 
This follows from Lemma \ref{016}. 
Indeed, $S=T_w$, where $w\in W$ is such that $w(\ep_i)=\ep_{i+1}$ for $i=1\ld n-1$, and $w(\ep_n)=-\ep_{1}$.
Note that $W$ acts on the set $\{\pm \ep_1\pm \cdots \pm\ep_n\}$ in the natural way as explained in Section 2.3,  and the $W$-orbit of $\om_i=\ep_1+\cdots +\ep_i$ consists of the elements $\pm\ep_{j_1}\pm \cdots\pm \ep_{j_i}$
($1\leq j_i\leq \cdots\leq j_i\leq n$).
Then it is clear that $w$ fixes none of these elements. It follows that
$w$ is not conjugate to any element of $W_j$ for $j=1\ld n$. $\Box$

\medskip
 Denote by $T^{\#}$ a maximal torus of $G$ whose cyclic decomposition is $T=T_1\times\cdots\times T_n$   and $|T_1|=\cdots =|T_n|=3$. The torus $T^{\#}$ is unique up to conjugacy as the Singer index of it must be $n$.
Then a canonical torus $D\subset \mathbf{T}$ consists of elements $\diag(t_1\ld t_n,t_n\up\ld t_1\up)$
with $t_i\in F, t_i^3=1$ for $i=1\ld n$. \itf $D$ is invariant under the Weyl group of $\mathbf{T}$. In addition,   $\ep_i(t)=t_i$ for $t\in T^{\#}$ and $i=1\ld n$.

\bl{ft2}
Let  $T=T_1\times \cdots \times T_k$ be a maximal torus of G of type $(\pm n_1\ld \pm n_k)$
with $n_1\geq n_2\geq \cdots \geq n_k$. 
Let $\lam=\om_1+\om_{2} $ and let $\chi_\lam$ be for the orbit character of $\lam$.

  $(1)$ $(\chi_{\lam}|_T,1_T)=0$ \ii $k=n$,
 $n_i=1$ for $i=1\ld n$  and $|T_1|=\cdots=|T_{n-1}|=3$.

 $(2)$   $1_T$ is not a constituent of $V_{\lam}|_T$ \ii $k=n$ and $T=T^{\#}$.

 $(3)$ Let $g\in G$ be a semisimple element. Then g has \ei $1$ on $V_{\lam}$.\el

Proof. Let $D$ be a canonical form of $T$, see Section 2.2. To simplify notation, we keep $T$
for $D$ and $t_i$ for $J_i$ for $i=1\ld k$. 

(1) Suppose first that $n_1>1$ and let $t\in T$.
 Then $\ep_1(t)=\zeta$, $\ep_2(t)=\zeta^2$, where $\zeta\in F$ with $|\zeta|=|T_1|$.  Then the $W$-orbit of $\lam=2\ep_1+\ep_2$ contains weight $\mu=2\ep_1-\ep_2$, and  $\mu(t)=1$.

 Next, suppose that $n_1=1$. Then $n_i=1$ and $|T_i|\in\{3,1\}$ for $i=1\ld n$. 
 If   $|T_{n-1}|=1$ then $|T_n|=1$ and $\ep_{n-1}(t)=\ep_{n}(t)=1$. As  the $W$-orbit of $\lam$ contains weight $\mu=2\ep_{n-1}+\ep_n$, we again get  $\mu(t)=1$.

Let  $|T_{n-1}|\neq 1$. The $W$-orbit of $\lam$ consists of weights $\pm 2\ep_i\pm \ep_j$  for any choice of $i,j$ with $i\neq j$. Then $ (\pm 2\ep_i\pm \ep_j)(T)\neq 1$.
Indeed, otherwise  $1=(\pm 2\ep_i\pm \ep_j)(T)= \{t_i^{\pm 2}t_j^{\pm 1}:  t_i\in T_i, t_j\in T_j\}$ for some choice of the signs, which is false.

(2)  By (1), we are left to inspect the cases where $k=n$,  $n_i=1$ for $i=1\ld n$  and $|T_1|=\cdots=|T_{n-1}|=3\geq |T_n|$.

 The only subdominant weights of $\lam$ are $\om_1,\om_3$ for $n>2$ and $\om_1$ for $n=2$. If $n>2$ then $\lam\in\Omega'$, and hence $\om_1,\om_3\in\Omega(V_\lam)$ by Lemma \ref{pr4}. If $n=2$ then $V_\lam=V_{\om_1}\otimes V_{\om_2}$, and hence  $\pm\ep_2+(\pm\ep_1\pm\ep_2)\in\Omega(V_\lam)$. So $\om_1=\ep_1\in\Omega(V_\lam)$. The $W$-orbit of $\om_1$ (respectively, of $\om_3$ if $n>2$) consists of weights
 $\pm \ep_1,...,\pm\ep_n$ (respectively, $\ep_{i_1}+\ep_{i_2}+\ep_{i_3}$ with $1\leq i_1<i_2< i_3\leq n$). If $|T_n|=1$ then $\ep_n(T)=1$. Let $|T_n|=3$. It is clear that $\ep_i(T)\neq 1_T$
 for every $i=1\ld n$, as well as $(\ep_{i_1}+\ep_{i_2}+\ep_{i_3})(T)\neq 1_T$ for $n>2$.  So  the result follows.

 (3)    If $g$ is contained in a maximal torus distinct from  $T^{\#}$ then $g$ has \ei 1 on $V_{\lam}$ by (2). This is obviously the case if  $g$ has \ei 1 on $V_{\om_1}$, the natural module for $G$.
  Otherwise, $g$ has no fixed point on $V_{\om_1}$.  Then $V_{\om_1}$ is a direct sum of $n$ two-dimensional
  $g$-stable subspaces, isomorphic to each other as $\FF_2\lan g\ran$-modules. In particular, $V_{\om_1}$ is a homogeneous $\FF_2\lan g\ran$-module. Then $C_G(g)\cong U_n(3)$, see for instance
  \cite[Lemma 6.6]{EZ}. Since $g$ is in the center of $U_n(3)$, this  lies in every maximal torus of $U_n(3)$. As every maximal torus of  $U_n(3)$ is a maximal torus of $G$ and $n>1$, there is a maximal torus $T\neq T^{\#}$ of $G$ with $g\in T$, so (3) follows from (2).  $\Box$

\bl{t37}
 Let $\om=\sum a_i\om_i\in\Omega^+$ be a   $2$-restricted weight. Suppose that  $a_n=1$ and $\delta(\om)\geq 2n$. Then $1_T$ is a constituent of $V_\om|_T$ for every maximal torus T of G.\el

Proof. If $\om$ is radical then $V_\om$ has weight 0 (Lemma \ref{t29}(1)), whence the result.
If $\om$ is not radical then, by Lemma \ref{t29}(2), $3\om_1\in\Omega(V_\om)$.
If $T=T^\#$ then  $t^3=1$ for every $t\in T^\#$, and hence $(3\om_1)(t)=\om_1(t^3)=1$, whence the result.

Let $T\neq T^\#$. By Lemma \ref{t29}(2),
 $\Omega(V_{\om_1+\om_2})\subseteq \Omega(V_\om)$. So the result follows from Lemma \ref{ft2}.  $\Box$

\bl{cc2}
 Let $\om=\sum a_i\om_i$ be a non-radical
$2$-restricted dominant weight and $a_n=0$. Suppose that $\om\neq \om_i$ for $i=1\ld n-1$ odd.
 
$(1)$ Let $T\neq T^{\#}$ be a maximal torus of $G$.
 Then $1_{T} $ is a constituent of $V_\om|_{T}$.

$(2)$ Let $g\in G$ be a semisimple element. Then g has \ei $1$ on $V_{\om}$.
\el

Proof. By Lemma \ref{sd4}, $\lam=\om_1+\om_2\preceq \om$.  In addition, all weights of $\Omega( V_{\lam})\subseteq\Omega(V_{\om})$  (Lemma \ref{pr4}). So the result follows from Lemma \ref{ft2}.  $\Box$

    \bl{t33}
     Let $ \om=\sum a_i\om_i$ be a $2$-restricted    dominant weight and $T=T^\#$.

    $(1)$ Suppose that  $a_n=0$.  Then $1_{T} $ is a constituent of $V_\om|_{T}$ 
    \ii either $0\preceq\om$ or $\sum a_i>2$.

    $(2)$   Suppose that  $a_n=1$. Then $1_{T} $ is a constituent of $V_\om|_{T}$
     \ii $\delta(\om)\geq 2n$.  \el

    Proof. Taking $T^\#$ to the canonical form (see Section 2.2), we can assume that $ t   =\diag(t_1\ld t_n,\\ t_n\up\ld t_1\up)$, where $t_1^3=\cdots = t_n^3=1$; then $\ep_i(t)=t_i$ for $i=1\ld n$.

    (1) The `only if' part. 
  Suppose the contrary. So $\om$ is not radical and   $\sum a_i\leq2$. Let $\mu=\sum a_i'\om_i$ be a dominant way    such that $\mu\prec \om.$  By Lemma \ref{sd3}(3), we have $\sum a_i'\leq 2$. Let  $\mu'$ be in the $W$-orbit of $\mu$. We claim that $\mu'(T^{\#})\neq1$ (getting a contradiction).    We can assume that $\mu=\om_i+\om_j$ $(i<j)$ as the cases with $\mu=\om_i$ and $\mu=2\om_i$ ($1\leq i\leq n-1$) are trivial. Then
     $\mu=2\ep_1+\cdots+2\ep_i+\ep_{i+1}+\cdots+\ep_j$ and $\mu'=2(\sum_{k\in K} \pm\ep_k)+(\sum_{l\in L} \pm\ep_l)$, where $K,L$ are disjoint subsets of $\{1\ld n\}$
  and $|K|=i,|L|=j-i$. Now the claim is obvious. 

  The `if' part. The case where $\om$ is radical is trivial. Suppose that $\om$ is not radical.
     By Lemma \ref{sd4}(3), $3\om_1$ is a subdominant weight of $\om$. Therefore, by Lemma \ref{pr4}(2),  $3\om_1$ is a weight of $V_\om$.    Obviously, $3\om_1(T^{\#})=1$, whence the result.

(2) Set $\om'=\om-\om_n$. Then we have $V_\om=V_{\om'}\otimes V_{\om_n}$.
So the weights of $V_\om$ are of shape $\nu+(\pm\ep_1\pm\cdots\pm\ep_n)$, where $\nu$ is a
  weight of $V_{\om'}$. Note that the $W$-orbit of $\nu$ contains a dominant weight
  and $w(T^{\#})=T^\#$. So we can assume that $\nu$ is dominant.

The `only if' part. Suppose  that $\delta(\om)<2n$. Then $\delta(\om')<n$.
By Lemma \ref{sd4}(6), $\delta(\nu)<n$ whenever   $\nu\prec\om'$ so $\om_n$ is not a subdominant weight of $\om'$.
Let $t\in T^\#$.  Then $((\nu+(\ep_1\pm\cdots\pm \ep_{n-1})) \pm\ep_n)(t)
=x\cdot (\ep_n(t_n)^{\pm 1}$, where $x=(\nu+(\ep_1\pm\cdots\pm\ep_{n-1}))(t)$. We can fix $t_1\ld t_{n-1}$ and
vary $t_n$ to be any $3$-root of 1. So we cannot  have $x\cdot (\ep_n(t_n))^{\pm 1}=1$  for arbitrary $t$.

The `if' part follows from Lemma \ref{t37}. $\Box$

\begin{theo}\label{ee3}
 Let $ \om=\sum a_i\om_i$ be a $2$-restricted  dominant weight. Then $1_{T}$ is
a constituent of $V_\om|_{T}$ for every maximal torus of G \ii one of the following holds:

$(1)$  $a_n=0$ and either $\sum a_i>2$ or $\delta(\om)=\sum a_ii$ is even;

$(2)$ $a_n=1$ and $\delta(\om)\geq 2n.$ 
\end{theo}

Proof.  The `if' part. In (1) if $\delta(\om)$ is even then $\om$ is radical,
and  $V_\om$ has weight 0 by Lemma \ref{pr4} (as $a_n=0$). This also holds in (2) by Lemma \ref{t29}.
So  the result follows in these cases. If $\sum a_i>2$ then the result follows from Lemma \ref{cc2}(1)  if $T\neq T^\#$ and from Lemma \ref{t33}(1) if $T=T^\#$.

The `only if' part.  Let $a_n=0$, and suppose  that $\sum a_ii$ is odd; then $\om$ is not radical. So $\sum a_i>2$ by Lemma \ref{t33}(1). If $a_n=1$ then the result follows  from Lemma \ref{t33}(2). $\Box$

\med{\it Proof of Theorem} \ref{ff3}.  By Lemma \ref{aa1}, we can assume that $|A|$ is odd. Therefore, it suffices to prove the result for maximal abelian groups of odd order. By Lemma \ref{ss1}, these are maximal tori of   $G$.
So the result follows from Theorem \ref{ee3}. $\Box$

\med Next we turn to proving Theorem \ref{si1}. For this we study the problem of the occurrence of the \ei 1 for individual elements $g\in G$ in \ir \reps $\rho$ that do not satisfy the condition (1), (2) of  Theorem \ref{ff3} for $A=\lan g \ran$. We start with the case where the \hw of $\rho$ is $\om_n$.

\def\se{semisimple }
\def\sel{semisimple element }
\def\selt{semisimple element. }
\def\mat{maximal torus }

Let $g\in G$ be a \sel   and $V$ the natural $\FF_2G$-module.  Set $A=\lan g \ran$, and let
$V=V_1\oplus\cdots \oplus V_k$ be a  decomposition of $V$ as a orthogonal sum of minimal non-degenerate $A$-stable subspaces, see Section 2.2.
Let $g_i$ be the restriction of   $g$ to $V_i$ for $i=1\ld k$, so we can write $g=\diag(g_1\ld g_k)$. Let $d_i=\dim V_i/2$. The decomposition $V=V_1\oplus\cdots \oplus V_k$ is unique up to reordering
of $V_1\ld V_k$ (if $|g|=3$ and $V_1$ is the sum of two totally isotropic $g$-stable subspaces then
$V_1$ is also the sum of two non-degenerate subspaces, so $V_1$ is not minimal.) Note that $g\in T$ for some 
\mat $T=T_1\times \cdots \times T_k$ of $G$ such that $TV_i=V_i$,  $g_i\in T_i$  for $i=1\ld k$, and $|T_i|=1$ if $g_i=1$. Then 
such a torus is below called {\it compatible with} $g$.

Define a graph $\Gamma(g)$  with vertices $1\ld k$, where the vertices $i,j\in\{1\ld k\}$ are linked \ii $(|g_i|,|g_j|)\neq 1$. A vertex $i$ is called {\it singular} if it is isolated and $|g_i|=2^{d_i}+1$;
 then $V_i$ is \ir $\FF_2A$-module. The set of  singular vertices is denoted by $\Gamma_0$ and
the number ${\rm Si}\,(g)$ of them is called the {\it Singer index of} $g$. So ${\rm Si}\,(g)=|\Gamma_0|$.
If  $i$ is singular then $\lan g_i\ran$ is a Hall subgroup of $A$ (that is, $(|g_i|,|A|/|g_i|)=1$).
If $T$ is a \mat compatible with $g$ then $|g_i|=|T_i|$ whenever $i$ is singular.
 
 \bl{e42}
 {\rm\cite[Proposition 4.12]{HZ09}} Let $\rho=\rho_{\om_n}$ and   let $T=T_1\times \cdots \times T_k$ be a maximal torus of G.
Set $m_i=\eta_i+1$. Then $\rho|_T=\otimes (\rho^{reg}_{T_i}+(-1)^{m_i} 1_{T_i})$, where $\rho^{reg}_{T_i}$ denotes the regular representation of $T_i$ and $1_{T_i}$ the trivial one. In particular, an \irr $\tau$ of $T$ occurs in $\rho|_T$ \ii $\tau(T_i)\neq 1_{T_i}$ whenever $\eta_i=-1$. \el

\begin{corol}\label{ch3}
Under notation of Lemma {\rm\ref{e42}} set $\rho_i'=\rho^{reg}_{T_i}+(-1)^{m_i} 1_{T_i}$.
Let $g\in T$ and $g=g_1\cdots g_k$ with $g_i\in T_i$ for $i\in\{1\ld k\}$.

 $(1)$  $\rho_i'(g_i)$ has less than $|g_i|$ distinct \eis \ii $\eta_i=-1$ and $|g_i|=|T_i|$.

 $(2)$ If $1$ is an \ei of $\rho_i'(g_i)$ then $\rho_i'(g_i)$ has $|g_i|$ distinct eigenvalues.

  $(3)$ the \eis of  $\rho_{\om_n}(g)$ are the products $e_1\cdots e_k$, where $e_i$ runs over the \eis of $\rho_i'(g_i)$ (or $\rho(g_i))$. Consequently, every primitive $|g|$-root of unity is an \ei of $\rho_{\om_n}(g)$.
\end{corol}

Proof. This follows straightforwardly from Lemma \ref{e42}. $\Box$

\bl{gg1}
 Let $g\in G$ be a \sel and let $\Gamma(g)$ be a graph defined above. Then the \f are equivalent:

$(1)$ $\rho_{\om_n}(g)$ has $|g|$ distinct eigenvalues;

$(2)$ $\rho_{\om_n}(g)$ has \ei $1;$

$(3)$ $\Gamma(g)$ has no singular vertex.
\el

 Proof. Set $\rho=\rho_{\om_n}$. For a connected component $J$ of  $\Gamma(g)$ set $g_J=\Pi_{i\in J}g_i$ and $\rho'_J=\otimes_{i\in J} \rho'_i$. We first prove

\med

$(*)$   $\rho'_J(g_J)$ has $|g_J|$ distinct eigenvalues unless $J=\{j\}$ is a singular vertex.

\med
Indeed,  if $|J|=1$ then $(*)$ follows from Corollary \ref{ch3}(1). Let  $|J|>1$.
Suppose first that $|J|=2$, and we can assume $J=\{1,2\}$.  By Lemma \ref{e42}, all non-trivial $|g_1|$-roots of unity $\lam_1\ld \lam_m$ are \eis of  $\rho'_1(g_1)$ (so $m=|g_1|-1$), and, by a similar reason, $\rho'_2(g_2)$ has at least 2 distinct \eis $\eta_1,\eta_2$, say, that are $|g_1|$-roots of unity. Then $\lam_i\eta_j$ $(1\leq i\leq m,j=1,2)$ yield all $|g_1|$-roots of unity, and they are \eis of $\rho'_1(g_1)\rho'_2(g_2)$. Similarly, all $|g_2|$-roots of unity are \eis of $\rho'_1(g_1)\rho'_2(g_2)$. As $\rho'_1(g_1)\rho'_2(g_2)=(\rho^{reg}_{T_1}(g_1)+(-1)^{m_i} 1_{T_1}(g_1))(\rho^{reg}_{T_2}(g_2)+(-1)^{m_i} 1_{T_2}(g_2))$, it follows that all other $|g_1g_2|$-roots of unity are \eis of $\rho'_1(g_1)\rho'_2(g_2)$, whence the claim. Let $|J|>2$. Then we argue similarly using induction on $|J|$.
 
$(3)\ra (1)$  follows from $(*)$ and Corollary \ref{ch3}(3), and $(1)\ra (2)$ is trivial.

$(2)\ra (3)$ Suppose the contrary, and let $i$ be a singular vertex. Then $|g_i|=|T_i|$ and $\eta_i=-1$.  \itf 1 is not an \ei of $\rho_i'(g_i)$, as well as of  $\rho(g)$ by Corollary \ref{ch3}(3) as  $|g_i|$ is coprime to $gg_i\up$.  $\Box$

\med
The \f lemma describes the set of \eis of $\rho_{\om_n}(g)$:

\bl{ff2}
   Let $g\in G$ be a semisimple element and $\Gamma _0$ the set of singular vertices on $\Gamma(g)$. For a $|g|$-root of unity  $\zeta\in \overline{\FF}_q$ let $|\zeta|$
denote the order of $\zeta$ in $\overline{\FF}_q^\times$. Then $\zeta$ is  an \ei of
$\rho_{\om_n}(g)$ \ii $(|\zeta|, |g_i|)>1$ for every $i\in \Gamma_0$. \el

 Proof. Let $T=T_1\times\cdots\times T_k$ be a \mat compatible with  $g$.
 Suppose $\zeta$ is an \ei of $\rho(g)$. By Corollary 3.9, $\zeta=e_1\cdots e_k$, where $e_i$ is an \ei of $\rho'_i(g_i)$ for $1\leq i\leq k$. So $e_i=\lam_i(g_i)$ for some \ir constituent $\lam_i$ of
$\rho'_i|_{T_i}$. Let $i\in \Gamma_0$.
Then $\lan g_i\ran=T_i$  and $\eta_i=-1$; therefore  1 is not an \ei of $\rho_i'(g_i)$ (Corollary 3.9). So  $e_i\neq 1$. Furthermore, $(|g_i|,|g_j|)=1$ for $j\neq i$, and hence $(|e_i|,|e _j|)=1$. Therefore, $|e_i|$ divides $|\zeta|$, as required.

Conversely, we have to show that $\zeta=e_1\cdots e_k$, where $e_i$ is an \ei of $\rho'_i(g_i)$ for $1\leq i\leq k$. Set $A=\lan g\ran$. As $|\zeta|$ divides  $|g|$ and $\lan g_i\ran$ is a Hall subgroup of  $A $ for $i\in \Gamma_0$, it follows that $\zeta$
is uniquely expressed as $\zeta_1\zeta_2$, where $|\zeta_2|$ divides $| g_i|$, and $(|\zeta_1|,| g_i|)=1$.
Then $\zeta_2\neq 1$ and we may set $e_i=\zeta_2$ (as all non-trivial $|g_i|$-roots of unity are \eis of
 $\rho_i'(g_i))$. Furthermore, 
 by reordering of $T_1\ld T_k$ we can assume that $\Gamma_0=\{1\ld r\}$, where
 $r$ is the Singer index of $g$. Set
  $D=\lan g_1\ld g_r\ran $; then $D$ is a Hall subgroup of $A$  and $D=T_1\times \cdots\times T_r$. So $A=D\times D_1$,
 where $(|D|,|D_1|)=1$.   Then $g=g_0h$, where $g_0=g_1\cdots  g_r\in D$ and $h=g_{r+1}\cdots g_k\in D_1$.
 The \eis of $\rho(g)$ are product of those
 of $\Pi_{i=1}^r \rho_i'(g_i)$ and of $\Pi_{i=r+1}^k \rho_i'(g_i)$. 
 It follows from  $(*)$ that  $\Pi_{i=r+1}^k \rho_i'(g_i)$ has $|h|$ distinct eigenvalues. In turn,
 the \eis of $\Pi_{i=1}^r \rho_i'(g_i)$ are exactly the products of non-trivial $|g_i|$-roots of unity for $i=1\ld r$. Whence the result. $\Box$

\med
Next we consider \reps of shape $\tau\otimes \rho_{\om_n}$, where $\tau\in\Irr {\mathbf G}$. We first observe the following:

  \bl{ef4}
    Let $\tau$ be a \rep of ${\mathbf G}=Sp_{2n}(\overline{\FF}_2)$ and     let $T$ 
 be  a maximal torus of G. Then $\rho |_T$ contains $1_T$ \ii $\tau$ has weight $\lam$ such that $\lam|_T=\om_n|_T$.
\el

Proof. Clearly, $1_T$ occurs in $(\tau\otimes \rho_{\om_n})|_T$ \ii there is an \ir constituent $\nu$, say, of $\rho_{\om_n}|_T$ such that $\nu\up$ occurs as an \ir constituent of $\tau|_{T}$.
 Let $\mu$ be a weight of $\tau$ such that $\nu=\mu(T)$. The Weyl group of $Sp_{2n}(F)$ contains $-\Id$, so $-\mu$ is a weight of $\tau$, and then $(-\mu)|_T\cong \nu\up$. Therefore, $\nu$ is a constituent of $\tau|_{T}$ too. As the weights of  $\rho_{\om_n}$ are exactly those  occurring in the $W$-orbit of $\om$, the lemma follows. $\Box$

\med
 Moreover,  if $\tau$ has weight $\om_n$ then $(\tau\otimes\rho_{\om_n})|_T> \rho^{reg}_T$ for any \mat $T$ of $G$. (We write
 $\si>\tau$ for \reps of $T$ if $\si=\tau+\nu$ for a proper \rep $\nu$ of $T$.) This essentially follows from Lemma \ref{e42} as $(\tau\otimes\rho_{\om_n})|_T> \otimes_i(\rho'_i\otimes \rho_i')$ and $(\rho'_i\otimes \rho_i')_{T_i}> \rho_{T_i}^{reg}$. 

\med
Similarly we have:

\bl{op2}
Let ${\mathbf G},\tau,\rho$ be as in Lemma {\rm \ref{ef4}}
 and $g\in G=Sp_{2n}(2)$ a semisimple element. Then $ \rho(g)$ has \ei $1$ \ii $\tau(g)$ and $\rho_{\om_n}(g)$ have a common eigenvalue. The latter holds \ii $\tau(g)$ has \ei $\zeta$ such that $(|\zeta|, |g_i|)>1$ for every $i\in \Gamma_0$.\el

Proof. The \eis of $(\tau\otimes \rho_{\om_n})(g)$ are products of those of $\tau(g)$ and $ \rho_{\om_n}(g)$. It is well known that all elements of $G$ are real, so if $\zeta$ is an \ei of $\tau(g)$ then so is $\zeta\up$. Whence the first statement of the lemma. The second one follows from Lemma \ref{ff2}. $\Box$

\med
So we are faced to decide when  $ \rho_{\om_n}(g)$ and
$\tau(g)$ have a common eigenvalue. This will be done in Section 4. We complete this section with a few useful observations.

\begin{corol}\label{e66}
 Let ${\mathbf G},\tau,\rho,g$ be as in Lemma {\rm \ref{op2}}.
Let h be the product of all $g_i$ with $i\in \Gamma_0$. Then $\rho(g)$ has \ei $1$ \ii $\rho(h)$ has \ei $1$.
\end{corol}

Proof. As $(|g_i|,|gg_i\up|)=1$  for $i\in \Gamma_0$, it follows that $h\in \lan g \ran$. Moreover,
$ \Gamma_0=\Gamma(h)$.   So the conditions in Lemma \ref{ff2} for $\rho(g)$ and $\rho(h)$ to have \ei 1 coincide, and the result follows. $\Box$

\begin{corol}\label{wwn}
Let ${\mathbf G},\tau,\rho,g$ be as in Lemma {\rm \ref{op2}}.
If $\tau$ has weight $\om_n$ then $\rho(g)$ has exactly $|g|$ distinct eigenvalues.
\end{corol}

Proof. Let $Q,P$ be the sets of all, all primitive $|g|$-roots of unity, respectively. By Corollary \ref{ch3}(3), every $e\in P$ is an \ei of $\rho_{\om_n}$. As $e_1e_2\in P$ whenever $e_1\in P$, $e_2\in Q\setminus P$, it follows that all non-primitive $|g|$-roots of unity are \eis of $\si(g)$.

As $|g|$ is odd, it follows that every $e\in P$ is a product $e_1e_2$ with $e_1,e_2\in P$.
Indeed, let $e_3\in Q\setminus P$. Then $e_3e,e_3\up e\in P$ and $(e_3e)(e_3\up e)=e^2\in P$.
There is a 2-power $k$, say, such that $e^{2k}= e$. Then $(e_3e)^{2k},(e_3\up e)^{2k}\in P$ and $(e_3e)^{2k}(e_3\up e)^{2k}=e$, as stated. $\Box$

\med
Remark. In fact we have shown that if $r$ is odd then every $r$-root of unity is a product of
two primitive $r$-roots of unity.

\medskip
{\it Proof of Theorem} \ref{th1}. If $a_n=0$ then the result follows from Lemma \ref{cc2}(2). If $a_n=1$ then this follows from   Lemma \ref{op2}, as $|\Gamma_0|\leq 1$ whenever $|g|$ is a prime power. (More precisely, if $|\Gamma_0|= 1$  then $\tau(g_i)\neq \Id$ (as $\tau$ is non-trivial by assumption), and hence $\tau(g_i)$ has an \ei $\zeta$ with $(|\zeta|,|g_i|)>1$ as required in Lemma \ref{op2}.) $\Box$

\section{Tensor products}

In this section $F$ is an \acf of arbitrary characteristic $p\geq 0$. Note that the results of Subsection 2.1
are independent from characteristic, so we are free to use them here.

\med
Let ${\mathbf G}=Sp_{2n}(F)$ and let ${\mathbf T}$ be the reference torus, that is, a maximal torus ${\mathbf T}=(t_1\ld t_n)$ $(t_1\ld t_n\in F^\times )$ such that $\ep_i({\mathbf T})=t_i$. Set ${\mathbf T}_i=\{(1\ld 1,t_i,1\ld 1): t_i\in F^\times \}$.

\bl{ob1}
 Let $\om\in\Omega^+$, $l=\delta(\om)$ and $k=\min(l,n)$. Let $V_\om$ be an
\ir ${\mathbf G}$-module with p-restricted \hw $\om$. Then there is a  p-restricted   dominant weight $\mu$ of $V_\om$ such that $\mu({\mathbf T}_i)\neq 1$ for $i=1\ld k$.\el

Proof. For uniformity set $\om_{n+1}=\om_1+\om_n$. By Lemma \ref{sd4}(6), $\om_l\preceq\om$ if $\l\leq n$, otherwise $\om_n\preceq\om$ or $\om_{n+1}\preceq\om$ if $l-n$ is even or odd, respectively.

Recall that $\om_i=\ep_1+\cdots+\ep_i$ for $i=1\ld n$ and $\om_{n+1}=2\ep_1+\cdots+\ep_n$. If $p\neq 2$
or $p=2$ and the $\om_n$-coefficient of $\om$ is 0, then, by Lemma  \ref{pr4},   $\om_l$, $\om_n$ or $\om_{n+1}$ is a weight of $\phi_\om$, provided, respectively, $l\leq n$, $l>n$ with $l-n$ even or $l>n$ with $l-n$ odd. We choose this weight for $\mu$. Then  $\mu({\mathbf T}_i)\neq 1$  for $i=1\ld k=\min(l,n)$, as stated.

Suppose that $p=2$ and the
$\om_n$-coefficient of $\om$ is 1. Then $V_\om=V_{\om'}\otimes V_{\om_n}$, where $\om'=\om-\om_n\in\Omega'$.
Then either $0\preceq\om'$ or $\om_1\preceq \om'$, so $V_{\om'}$ has weight 0 or $\om_1$ by  Lemma
 \ref{pr4}. Then $V_\om$ has weight $\mu=\om_n=\ep_1+...+\ep_n$ or $\mu=\om_1+\om_n=2\ep_1+...+\ep_n$.
So again $\mu({\mathbf T}_i)\neq 1$ for $i=1\ld n$.  $\Box  $

\begin{theo}\label{th7}
Let $\mathbf{ G}=Sp_{2n}(F)$ and let  
$n_1\ld n_k$ be positive integers with $ n_1+\cdots+n_k\leq n$. Let $\mathbf{ H}={\mathbf H}_1\times\cdots\times {\mathbf H}_k$, where ${\mathbf H}_i\cong Sp_{2n_i}(F)$ for $i=1\ld k$.
Let $\phi=\phi_\om$ be an  \irr of $\mathbf{ G}$ with p-restricted \hw $\om=\sum a_i\om_i$. Then the \f are equivalent:

$(1)$ $\rho({\mathbf H}_i)=1$ for some $i\in\{1\ld k\}$ whenever $\rho$ is a composition factor   of $\phi|_\mathbf{ H};$

$(2)$ $\delta(\om)=\sum a_ii<k$.

\noindent In addition, the result is valid for $\tilde \phi_\om$, the \rep afforded by the Weyl module with \hw $\om$.\end{theo}

Proof. It suffices to prove the theorem for the case where $n_1=\cdots =n_k=1$. Indeed, we can choose a subgroup $X_i=Sp_2(F)$ in every $\mathbf{ H}_i$, and set $X=X_1\times \cdots\times X_k$. Then $\rho({\mathbf H}_i)=1$ \ii $\rho(X_i)=1$, whence the claim.

$(2)\ra (1)$
Let $ {\mathbf T}'$ be a maximal torus of $\mathbf{ H}={\mathbf H}_1\times\cdots \times {\mathbf H}_k$ (where now ${\mathbf H}_1\cong \cdots \cong {\mathbf H}_k\cong Sp_2(F)$). Let ${\mathbf T}_i'$ be a maximal torus of ${\mathbf H}_i$, so ${\mathbf T}'={\mathbf T}_1'\times\cdots \times  {\mathbf T}_k'$. We can assume that ${\mathbf T}'\subset {\mathbf T}$ and ${\mathbf T}'=(t'_1\ld t'_k)$ in the sense that
$\ep_i((t'_1\ld t'_k))=t'_i$ for $i\leq k$ and $\ep_i((t'_1\ld t'_k))=1$ for $i> k$. Let $\si$ be a weight of $\phi$, and write $\si=\sum_{i=1}^n c_i\ep_i$ for some integers $c_i$.  Then $\si|_{{\mathbf T}'}=(c_1\ep_1\ld c_k\ep_k)$, that is, $\ep_i$ as a function on ${\mathbf T}_i'$ is  in fact the fundamental weight of ${\mathbf H}_i$ for $i=1\ld k$. If $\si$ is such that
 $\si|_{{\mathbf T}'}$ is a highest weight of an \ir constituent $\rho$ of $\phi|_\mathbf{ H}$ then $c_i\geq 0$ for $i=1\ld k$. Let $\nu=\sum d_i\ep_i$ be a dominant weight in the $W$-orbit $W\si$.
 As $\nu$ is a weight of $\phi$, we have $\nu\leq \om$, and then $\delta(\nu)\leq \delta(\om)<k$ by Lemma \ref{sd3}(3). Note that  $\delta(\nu)=d_1+\cdots +d_n$, so $d_1+\cdots +d_n<k$.
The action of $W$ on the weights preserves the set $\pm \ep_1\ld \pm\ep_n$; it follows that the set
$(|c_1|\ld |c_n|)$ coincides with $(d_1\ld d_n)$ up to reordering,   where
 $|c_i|$ means the absolute value  of $c_i$. So  $\sum_i |c_i|\cdot i<k$. As $c_1\ld c_k$ are non-negative, we have  $\sum_{i=1}^k c_i i<k$, and hence at least one of $c_1\ld c_k$ equals 0. This means that $\si({\mathbf H}_i)=1$ for this $i$,  and hence ${\mathbf H}_i$ is in the kernel of $\rho$, as required.

$(1)\ra (2)$ follows  from Lemma \ref{ob1} as ${\mathbf T}_i$ is a maximal torus of ${\mathbf H}_i$.

For additional statement, $(1)\ra(2)$ follows as $\phi_\om$ is a composition factor of  $\tilde \phi_\om$.
As $(2)\ra(1)$ is true for $F=\CC$, and the weights of $\tilde \phi$ are the same as those of the
\irr of $Sp_{2n}(\CC)$ with \hw $\om$, the above reasoning remains valid for $\tilde \phi_\om$. $\Box$

\begin{corol}\label{y10}
In notation of Theorem $\ref{th7}$ let $p>0$ and $\tau=\phi_1\otimes\cdots\otimes\phi_m$, where $\phi_j$ is an \irr of ${\mathbf G}$ with $p$-restricted \hw $\mu_j$, $j=1\ld m$. Let $l=\delta(\mu_1)+\cdots +\delta(\mu_m)$. Suppose that $l<n$. Let $\rho$ be an \ir constituent of $\tau|_{{\mathbf H}}$ and $J=\{i\in \{1\ld k\}: \rho({\mathbf H}_i)\neq 1\}$.
Then $|J|\leq l$. 
\end{corol}

Proof. Note that $\tau|_{{\mathbf H}}=\otimes \phi_j|_{{\mathbf H}}$, so $\rho=\otimes\si_j$, where
$\si_j$ is an \ir constituent of $\phi_j|_{{\mathbf H}}$. In turn, $\si_j=\nu_1\otimes\cdots\otimes \nu_k$,
where $\nu_i\in\Irr {\mathbf H}_i$ for $i=1\ld k$. By Theorem \ref{th7}, this product has at most
$\delta(\mu_j)$ non-trivial multiples. (Indeed, let $l_j=\delta(\mu_j)$; if the claim is false then $\si_j$ has at least $l_j+1$
non-trivial multiples, so $k>l_j$. We can assume that these are $\nu_1\ld \nu_{l_j+1}$. Set ${\mathbf H}'={\mathbf H}_1\times\cdots \times{\mathbf H}_{l_j+1}$ and $\rho'=\nu_1\cdots \nu_{l_j+1}$. By Theorem \ref{th7} applied to ${\mathbf H}'$ in place of ${\mathbf H}$ and $\rho'$ in place of $\rho$, there is $i\in \{1\ld l_{j+1}\}$ such that $\si'({\mathbf H}_i)=1$, which is a contradiction.)  Therefore, $\rho$ has at most $l=\sum _j\delta(\mu_j)$ non-trivial multiples. $\Box$

\med Recall that the Singer index of a maximal torus $T=T_1\times \cdots \times T_k$ of $G$ is the number of $i\in\{1\ld k\}$ with $\eta_i=-1$.

\begin{theo}\label{s10}
Let $G=Sp_{2n}(2)$ and let $T=T_1\times\cdots\times T_k$ be a maximal torus of G of Singer index $m$. Let $\om=(a_1\ld a_{n-1},a_n)$ with $a_n=1$
be a $2$-restricted weight and $\phi_\om$ be an \irr of G with \hw $\om$. Then $1_T$ is a constituent of $\phi_\om|_T$ \ii $\sum a_ii\geq n+m$.\end{theo}

Proof. Set $\om'=\om-\om_n$. By Lemma \ref{ef4}, $1_T$ is a constituent of $\phi|_{T}$ \ii $\phi_{\om_n}|_T$ and $\phi_{\om'}|_T$ have a common \ir constituent. By Lemma \ref{e42}, an \irr $\tau$ of $T$ occurs in $\phi_{\om_n}|_T$ \ii $\tau(T_i)\neq 1 $ whenever   $\eta_i(T_i)=-1$ $(1\leq i\leq k)$. We show that $\tau$ with this property is a constituent of $\phi_{\om'}|_T$ \ii $\delta(\om)\geq \sum _ia_ii\geq m+n$, equivalently, $\delta(\om')\geq m$. %

By reordering of $T_1\ld T_k$ we can assume that $\eta_i=-1$ for $i=1\ld m$. One easily observes that there is a subgroup ${\mathbf H}={\mathbf H}_1\times\cdots\times  {\mathbf H}_m$ of $Sp_{2n}(\overline{\FF}_2)$ such that $ {\mathbf H}_i=Sp_{2n_i}(\overline{\FF}_2)$ for $i=1\ld m$, and $T_i$ is a maximal torus of $Sp_{2n_i}(2)$. Let
$\overline{\phi}_{\om'}$ be a \rep of ${\mathbf G}=Sp_{2n}(\overline{\FF}_2)$ with \hw $\om'$. Then  $\overline{\phi}_{\om'}|_G=\phi_{\om'}$.

Suppose that $\delta(\om')\geq m$. By Theorem \ref{th7}, there is a composition factor of $\overline{\phi}_{\om'}|_{\mathbf H}$ which is non-trivial on each ${\mathbf H}_i$ for $i=1\ld k$. Then $\overline{\phi}_{\om'}|_{T_i}\neq 1$ for $i=1\ld m$ 
as such $T_i\neq 1$.   (Note that $\eta_i=- 1$ implies $|T_i|\neq 1$.) 
So $1_T$ is a constituent of $\phi_\om|_{T}$ by the above. 

Conversely, if $\phi_{\om'}|_{T}$ has a constituent $\tau$, say, non-trivial on each $T_i$ with $1\leq i\leq m$ then $\overline{\phi}_{\om'}|_{\mathbf H}$ has an \ir constituent that is non-trivial on ${\mathbf H}_i$
for $i=1\ld m$.   Then  $\sum _{i=1}^{n-1}a_ii>m$ by Theorem \ref{th7}, and the result follows. $\Box$

\begin{theo}\label{fr1}
 Let $g\in G $ be a semisimple element  and let $\rho=\rho_{\om'}\otimes \rho_{\om_n}$ be an \irr of G. Then $1$ is an \ei of
$\rho(g)$   \ii $\delta(\om)\geq n+ {\rm Si}(g)$.
\end{theo}

Proof.  Let $h$ be the product of $g_i$ with $i\in \Gamma_0$. Then $ {\rm Si}(h) ={\rm Si}(g)$.
 By Corollary \ref{e66}, $\rho(g)$ has \ei 1 \ii $\rho(h)$ has \ei 1.
 So we can assume that $g=h$. Then $T=\lan g\ran$ as in this case $(|g_i|,|g_j|)=1$ for all $i\in\Gamma(g)$
 and $g_i=1$ if $i\notin \Gamma_0$.  Therefore, the Singer index of $g$ equals that of $T$.
Then the result follows from Theorem \ref{s10}.  $\Box$

\med{\it Proof of Theorem} \ref{si1}. If $a_n=0$ then the result follows from  Lemma \ref{cc2}
(the "if" part) and Lemma \ref{op1} (the "only if" part).

Let $a_n=1$. Then ${\rm Si}(g)\leq {\rm Si}(n)$, and the equality holds for some $g\in G$. So the result follows from Theorem \ref{fr1}. $\Box$

\med
Observe that there are no elements whose Singer index  is greater than ${\rm Si}(n)$, and there are semisimple elements $g\in G$ whose  Singer index equals   ${\rm Si}(n)$.
We refer to such elements as those of {\it maximal Singer index}.

\begin{propo}\label{p99}
Let $G=Sp_{2n}(2)$ and let $t\in G$ be an element of maximal Singer index. Let $\rho$ be an \irr of G with  \hw $\om=(a_1\ld a_n)$, $a_n=1$. If $1$ is an \ei of $\rho(t)$ then $1$ is an \ei of $\rho(g)$ for every $g\in G$.\end{propo}

Proof. This follows from Theorem \ref{fr1}.  $\Box$

\begin{propo}\label{p88}
Let $g\in G=Sp_{2n}(2)$ and ${\mathbf G}=Sp_{2n}(\overline{\FF}_2)$. Suppose that $\rho_{\om_i}(g)$ has \ei $1$ for  $i=1, n$. Then $\rho(g)$ has \ei $1$  for every $\rho\in\Irr {\mathbf G}$.
\end{propo}

Proof. Let $\om=(a_1\ld a_n)$ be the \hw of $\rho$. Suppose first that $\om$ is 2-restricted. If $a_n=0$ and $\om$ is radical then $\rho$ has weight 0 by Lemma \ref{t29}(1), and the claim is trivial. If $\om$ is not radical then $\rho$ has weight $\om_1$ and $\Omega(\rho_{\om_1})\subset \Omega(\rho)$, whence the result.

Let $a_n=1$.  By assumption, $\rho_{\om_n}(g)$ has \ei 1, so  $\rho_{\om_n}(g)$ has $|g|$ distinct eigenvalues by Lemma \ref{gg1}. Let $\om'=\om-\om_n$. Then $\rho_{\om}=\rho_{\om'}\otimes \rho_{\om_n}$, so
$\rho_{\om}(g)$ has $|g|$ distinct eigenvalues, whence the result in this case.

Finally suppose that $\rho$ is not 2-restricted. Then $\rho|_G$ is a tensor product of \ir \reps of $G$, and we conclude similarly. $\Box$

\bl{tp1}\mar{tp1}
Let $V=V_{\om_i}\otimes V_{\om_j}$ for $1\leq i< j\leq n$.
Let $g\in G$ be a semisimple element.  If $j<n$ or
$j=n$ and $i\geq {\rm Si}(g)$ then g has \ei $1$ on V.
\el

Proof.  It is well known that $V$ has a composition factor $V_{\om_i+\om_j}$. If $j<n$ then the result follows by Theorem \ref{si1}(1). If $j=n$ then $V_{\om_i}\otimes V_{\om_n}\cong V_{\om_i+\om_n}$, and  the result follows by Theorem \ref{si1}(2).
  $\Box$

\med

Note that ${\rm Si}(g)= |\Gamma_0(g)|\leq {\rm Si}(n)$ and ${\rm Si}(n)<n$. In addition,  ${\rm Si}(3)=2={\rm Si}(4)={\rm Si}(5)={\rm Si}(6)$, ${\rm Si}(7)=3={\rm Si}(8)={\rm Si}(9)={\rm Si}(10)={\rm Si}(11)$, ${\rm Si}(12)=4$.

\begin{propo}\label{p49} 
Let $G=Sp_{2n}(2)$, ${\mathbf G}=Sp_{2n}(\overline{\FF}_2)$ and let $\rho$ be an \irr of ${\mathbf G}$ with \hw  $\om=\sum a_i\om_i$. 
Let $g\in G$ be a semisimple elt. Then one of the \f holds:

$(1)$ $1$ is an \ei of $\rho_\om(g);$

$(2)$ $\om= 2^i\om_j$ for  $j$ odd or $j=n$; 

$(3)$  $\om=\om'+2^k\om_n$ with $0\neq \om'\in \Omega'$ and $k\geq 0$. In addition, if
$\om'=\sum 2^i\nu_i$, where $\nu_i$'s are $2$-restricted, then   $\delta (\sum \nu_i)<{\rm Si}(g)$.
 \end{propo}

Proof. Suppose that (1),(2) do not hold. Let $\rho=\otimes_{i\geq 0} \rho_{2^i\mu_i}$, where $\mu_i$ are 2-restricted dominant weights.
Note that  $\rho_{2^i\mu_i}|_G\cong \rho_{\mu_i}|_G$. So $\rho(g)\cong (\otimes\rho_{\mu_i})(g)$.
Set $\tau=\otimes\rho_{\mu_i}$. If $\mu_i\notin \Omega'$ then $\mu=\nu_i+\om_n$, where $\nu_i\in\Omega'$, and  then $\rho_{\mu_i}=\rho_{\nu_i}\otimes\rho_{\om_n}$.
So  $\tau=\tau' \otimes \si$, where $\tau'=\otimes\rho_{\nu_i}$ with $\nu_i\in\Omega'$ and $\si$ is either trivial or the tensor product of $k>0$ copies of $\rho_{\om_n}$; the former is expressed as $k=0$.
As (1) does not holds,   $\tau$ does not have weight 0.

If $k>1$ then $(\rho_{\om_n}\otimes\rho_{\om_n})(g)$ has $|g|$ distinct \eis (Corollary \ref{wwn}), and hence 
so is $\rho(g)$. This remains true if $k=1$ and $\om_n$ is a weight of the product of any multiples of
 $\otimes\rho_{\nu_i}$ (as $\Omega(\rho_{\om_n})=W\om_n$). So, to prove (3), we can assume that $k\leq 1$.

Let $k=1$. So $\tau=\tau' \otimes \rho_{\om_n}$, where $\tau'=\otimes\rho_{\nu_i}$, and
we can assume that $\om_n$ is not a weight of any product of $\rho_{\nu_i}$;
 in particular, $\delta(\nu_i)< n$ by Lemmas \ref{sd4}(6) and \ref{pr4} (as $\nu_i\in\Omega'$). Set
$l_i=\delta(\nu_i)$, so $l_i<n$. Again by Lemma \ref{sd4}(6), $\rho_{\nu_i}$ has weight $\om_{l_i}$, and hence, by Lemma \ref{pr4}, the weights of  $\rho_{\om_{l_i}}$  are also weights of $\rho_{\nu_i}$. Note that $\delta(\nu_i)=\delta(\om_{l_i})$.

Observe that  $\om_j,\om_{j-2},\om_{j-4}\ldots$ are weights of $\rho_{\om_j}$ for $j<n$. \itf  $\rho_{\om_{i}}\otimes  \rho_{\om_{j}}$ for $i,j<n$ has all weights $\om_{r}$ for $r\leq \min\{n,i+j\}$ and $i+j-r$ even. (This is clear if $i\neq j$ and for $i=j$ with $2i>n$. Let $2i\leq n$. Then  $\rho_{\om_{i}}\otimes  \rho_{\om_{i}}$ has weights $2\om_i$ and $\om_{i-1}+\om_{i+1}$ as well as $\om_{2i}$.)

Suppose that $\sum l_i\geq n$, $k=1$.  Then either $\om_n$ or $\om_{n-1}$ is a weight of  $\tau'$. In the former case
  1 is observed above to be an \ei of $\tau(g)$. In the latter case $\om_{n-3},\om_{n-5},\ldots$
  are weights of $\tau'$;
 this implies  $\Omega(\rho_{\om_{n-1}})\subseteq \Omega(\tau')$.
  So it suffices to deal with the case
 where $\tau'=\rho_{\om_{n-1}}$, $\tau=\rho_{\om_{n-1}}\otimes\rho_{\om_{n}}$ (ignoring that $\delta(\om_{n-1})<n$). So $\delta(\tau)=2n-1$.
By Theorem \ref{fr1}, if $\delta(\tau)=2n-1\geq n +{\rm Si}(g)$ then 1 is an \ei of $\tau$. However,
${\rm Si}(g)\leq {\rm Si}(n)\leq n-1$ for $n>1$. So we are done.

Next assume that  $k=1$, $\sum l_i<n$, 
and let $l=\sum l_i$.  The above observation on weights of $\rho_{\om_{i}}\otimes  \rho_{\om_{j}}$ for $i,j<n$ implies $\om_l$ to be a weight of $\tau'$ as well as all $\om_{l-2},\om_{l-4}\ldots$. This means  $\Omega(V_{\om_l})\subseteq \Omega(\tau')$. By Lemma \ref{tp1}, if $l\geq {\rm Si}(g)$ then  1 is an \ei of $(V_{\om_l}\otimes V_{\om_n})(g)$, and hence
$\tau(g)$, whereas we have assumed that (1) does not hold. Therefore, $l< {\rm Si}(g)$, whence (3) in this case.

Let $k=0$ so $\tau=\tau'$ and $\mu_i=\nu_i$. By reordering of the multiples of $\tau$ we can assume that
$\pm\om_1$ are weights of $\rho_{\mu_1}$ (as $\tau$ does not have weight 0).
If $\tau=\rho_{\mu_1} $ then the result follows from Theorem \ref{th1}, so we assume that $\si:=\Pi_{i>1}\rho_{\mu_i}$ is not trivial. By Lemmas \ref{sd4} and \ref{pr4}, either $\pm \om_1$ or $0,\om_2$ are weights of every multiple. The tensor product of two multiples with the former property has weight 0. \itf that $0,\om_2\in\Omega(\si)$, and hence  $\om_1+\om_2 ,\om_1\in\Omega(\tau)$. If $n>2$ then $\ep_2+\ep_3\in\Omega(\tau)$, and hence $\om_3\in\Omega(\tau)$.
So  $\Omega(\rho_{\om_1+\om_2})\subseteq \Omega(\tau)$. Then 1 is an \ei of $\tau$ by Lemma \ref{ft2}(3).   $\Box$

{\it Proof of Theorem} 1.8. The "if" part follows from Lemma 3.2. The "only if" part follows from Proposition 4.9.
Indeed, we only have to show that (3) of Proposition 4.9 does not hold. For this,
in notation of 4.9 we have ${\rm Si}(g)=1$ (as defined prior to Lemma 3.8), and $\om'\neq 0$. Then $\sum \nu_i\neq 0$,
and hence $\delta(\sum \nu_i)\geq 1$. It follows that (3) fails, and we are done. $\Box$  


  \section{An application}

 \subsection{Remarks on Weyl modules}

For a definition of Weyl module we refer to \cite[\S 3.1]{Hum}. Here we mainly use the fact that, for a simple \ag $\mathbf{X}$, a Weyl module   is an  $\mathbf{X}$-module (in general, reducible)
whose weight set coincides  with the weight set  of some \irr of the Lie algebra of $\mathbf{X}$, equivalently, of the simple \ag over the complex numbers of the same type as that of $\mathbf{X}$
 (regarding the multiplicities of weights). (The word "coincides" here means the coincidence of the respective sets of strings $(a_1\ld a_n)$.) By \cite[Ch. VIII, Proposition 5]{Bo8}, this implies

\bl{gw1}
Let $\mathbf{X}$ be a simple algebraic group,   $\tilde{ V}_\om$   a Weyl module of $\mathbf{X}$ with \hw   $\om$ and $\mu\prec \om$  a dominant weight of  $\mathbf{X}$. Then $\mu$ is a weight of  $\tilde V_\om$.
\el

Let $\mathbf{G}=Sp_{2n}(F)$ and the characteristic of $F$ equal 2. Then the weight set  of  $\tilde V_\om$, $\om=\sum a_i\om_i$,
in contrast with \ir $\mathbf{G}$-module $V_{\om}$, has no singularity when $a_n\neq 0$.   In particular, we have the following version of Lemma \ref{cc2}(2)   for Weyl modules:

\bl{gc4}
Let $\om$ be a dominant weight of $\mathbf{G}$ and $\tilde{V}=\tilde{V}_\om$ a Weyl module of  highest weight $\om$. Let $g\in G$ be a  semisimple element. Then g has \ei $1$ on $\tilde{V}$ unless possibly $\om=\om_i$ for i odd.\el

Proof. If $\om$ is radical then the zero weight is a weight of $\tilde{V}$,  and the claim follows. Suppose $\om$ is not radical.  Then $\om_1\preceq \om$. As $\om\neq \om_i$, we have $ \om_1+\om_2\preceq \om$ by Lemma \ref{sd4}(5). By Lemma \ref{pr4}(2), 
 $\om_1+\om_2\in\Omega(\tilde{V})$ and $\Omega(V_{\om_1+\om_2})\subseteq \Omega(\tilde{V})$. 
 So Lemma \ref{ft2}(3) yields the result.
 
 \subsection{$SL_{2n}(F)$-modules viewed as $Sp_{2n}(F)$-modules}

Let $F$ be an \acf and $\mathbf{H}=SL_{2n}(F)$, $\mathbf{G}=Sp_{2n}(F)$.  Let $\mathbf{T}$, $\mathbf{T}'$ be maximal  tori of  $\mathbf{H}$, $\mathbf{G}$, respectively.
Under a Witt basis  $e_1\ld e_{2n}$ of the natural $F\mathbf{G}$-module, we assume that
$\mathbf{T}$  consists of matrices $t=\diag(t_1\ld t_{2n})$ with $t_1\ld t_{2n}\in F^\times $ and $t_1\cdots  t_{2n}=1$, and $\mathbf{T}'$ consists of matrices $t'=\diag(t_1\ld t_n,t_n\up\ld  t_1\up$)   ($t_1\ld t_n\in F^\times )$.
So $\mathbf{T}'\subset \mathbf{T}$. The Bourbaki weights $\ep_1\ld \ep_{2n}$  for $\mathbf{H}$
are defined by $\ep_i(t)=t_i$ ($i=1\ld 2n$). We denote by $\ep'_1\ld \ep'_{n}$   the Bourbaki weights for $\mathbf{G}$, which are defined by $\ep'_i(t')=t_i$ for $i=1\ld n$.

Then $\ep_i|_{\mathbf{T}'}=\ep'_{i}$ if $i\leq n$ and $\ep_i|_{\mathbf{T}'}=-\ep'_{2n+1-i}$ if $i>n$ (or $\ep_{2n-i}|_{\mathbf{T}'}=-\ep'_{i+1}$ for $1\ld n-1$).
We denote $\om_1 \ld \om_{2n-1}$ the fundamental weights for $\mathbf{H}$
and  $\om'_1\ld \om'_{n}$  the  fundamental weights for $\mathbf{G}$.

\med
In Lemmas \ref{w0}, \ref{ww1} and \ref{ww0}  $F$ is a field  of arbitrary characteristic $p>0$.

 \bl{w0}
  Let p be the characteristic of F. Let V be an \ir $\mathbf{H}$-module with $p$-restricted highest weight, and $\mu$ a dominant weight of $V|_{ \mathbf{G} }$ (not necessarily $p$-restricted).
Then all weights of  $\tilde V_\mu$ occur as weights of $V|_{\mathbf{H}}$.
 \el

Proof. By Lemma \ref{gw1} and \cite{Su}, the set of weights of $V$ coincides with that of an \ir $\tilde{\mathbf{H}}$-module $\tilde V$ whose \hw is the same as that of $V$. So $\mu$ is a weight of
$\tilde{V}|_{\mathbf{G}}$. Let $L_{\mathbf{H}}$ and $L_{\mathbf{G}}$ be the Lie algebra over the complex numbers of the same type as $\mathbf{H}$ and $\mathbf{G}$, respectively. (So $L_{\mathbf{H}}$ is of type $A_{2n-1}$ and $L_{\mathbf{G}}$ is of type $C_{n}$.) Then $\mu$ as a weight of $L_{\mathbf{G}}$ is a dominant weight of some \irr $\Lambda$, say,  of $L_{\mathbf{G}}$, and all weights of $\Lambda$ are weights of $\tilde{V}|_{\mathbf{G}}$. By \cite[Ch. VIII, Prop. 5]{Bo8}, all weights of
\irr of $L_{\mathbf{G}}$ with \hw $\mu$ are weights of $\Lambda$. Therefore, all weights of $\tilde V_\mu$
are weights of $\tilde{V}|_{\mathbf{G}}$.

 \bl{ww1}
 Let $\om=\sum a_i\om_i$ be a dominant weight of $\mathbf{H}$. Then $\om|_{\mathbf{T}'}=\sum b_i\om_i'$, where $b_i=a_i+a_{2n-i}$ for $i=1\ld n-1$ and $b_n=a_n$.
 \el

 Proof. It suffices to prove the lemma for $\om=\om_i$. If $i\leq n$ then $\om_i|_{T'}=(\ep_1+\cdots +\ep_i)|_{T'}=\ep'_1+\cdots +\ep'_i=\om_i'$. Let $i>n$. Then $\om_i|_{T'}=(\ep_1+\cdots +\ep_i)|_{T'}=\ep'_1+\cdots +\ep'_n-\ep'_{n-1}-\cdots -\ep_{2n-i+1}=\ep_1'+\cdots+\ep_{2n-i}'=\om_{2n-i}'$. $\Box$

\bl{ww0}
 Let V be an \ir $\mathbf{H}$-module with highest weight $\om_k$ for $1\leq k\leq 2n-1$,
 and let $V'$ be a composition factor of $V|_{\mathbf{G}}$. Then the highest weight of $V'$ is $\om_j'$
 for some $j$ with $1\leq j\leq \min(k,2n-k)$ and $k-j$ even. \el

 Proof. We can assume $k\leq n$ by replacing $V$ by the dual of it. The weights of $V$ form the $W$-orbit of $\om_k=\ep_1+\cdots+\ep_k$, where $W$ is the Weyl group of $\mathbf{H}$. So if $\mu$ is a weight of $V$ then $\mu=\ep_{i_1}+\cdots+\ep_{i_k}$, where $1\leq i_1\leq ...\leq  i_k\leq  2n$. Let $l$ be an integer such that $i_{l-1}\leq n<i_l$ (it exists unless $i_k\leq n$ or $i_1>n$).
 Then $\mu|_{\mathbf{T}'}=\ep'_{i_1}+\cdots+\ep'_{i_{l-1}}-\ep'_{2n+1-i_l}-\cdots-\ep_{2n+1-i_k}$ (with an obvious refining for the exceptional cases). After canceling the terms occurring here with opposite signs, one obtains
a similar expression $\mu'=\sum b_r \ep'_r$ with $-1\leq b_r\leq 1$ for $1\leq r\leq k$. Moreover,
 $k-s$  is even, where $s$ is the number of non-zero coefficients $b_r$. This is a weight of some \ir constituent of $V|_{\mathbf{G}}$.  Recall that $W'$, the Weyl group of $\mathbf{G}$, acts transitively on
every set $\pm\ep'_{m_1}\pm\cdots\pm\ep'_{m_s}$ for any fixed $s$ and integers  $1\leq m_1\leq ...\leq  m_s\leq  n$.
Then the dominant weight in the orbit $W'\mu'$ is $\ep'_{1}+\cdots+\ep'_{s}$. We conclude
that the highest weight of $V'$ must be $\ep'_1+\cdots+\ep'_j=\om_j $ for some $j$, as stated.


\begin{corol}\label{no1}
Let $\om=\sum a_i\om_i$ be a $2$-restricted dominant  weight and $\om'=\om|_{\mathbf{T}'}=\sum a'_i\om_i'$.

$(1)$ $a_i'\leq 2$ for $i=1\ld n-1$ and $a_n'\leq 1;$

$(2)$ if $\om$ is not radical then $\om'\neq 2\om_j'$ for any $j=1\ld n;$

$(3)$ if $\om'=\om'_j$ for some j then $\om=\om_j$ or $\om_{2n-j}$.
\end{corol}

Proof.
(1) follows from Lemma \ref{ww1}. (2) If $j=n$ then the claim is a part of (1). Let $j<n$. Then $\om'= 2\om_j'$ yields $a_j+a_{2n-j}=2$ and $a_i+a_{2n-i}=0$ for $i\neq j$ by Lemma \ref{ww1}. This means
that $\om=\om_j+\om_{2n-j}$, which is a radical weight. This is a contradiction.

(3)   If $j<n$ then $a_j+a_{2n-j}=1$ and $a_i+a_{2n-i}=0$ for $i\neq j$, $a_n=0$.
 If  $j=n$ then $a_n=1$ and $a_i+a_{2n-i}=0$ for $i<n$. So the claim follows. $\Box$

 \bl{do1}\mar{do1} Let $V$ be an \ir $\mathbf{H}$-module with $2$-restricted highest weight $\om$.
 If  $\om\neq \om_i$ for any $i\in\{1\ld 2n-1\}$    then $V|_{\mathbf{G}}$ has either weight $0$ or all weights of $V_{\om'_1+\om_2'}.$   \el

 Proof.   Set $\om'=\om|_{\mathbf{T}'}$. It follows from \cite[Theorem 39]{St} that $V|_{\mathbf{G}}$ has a composition factor $U$, say, of \hw $\om'$. (Indeed, let $v^+$ be a vector of weight $\om$. Then $v^+$
 is stabilized by the unipotent radical of certain Borel subgroup $B$ of $\mathbf{H}$, and hence
 by the unipotent radical of a Borel subgroup $B'$ of $\mathbf{G}$ contained in $B$. Set $V'=  \mathbf{G}v^+$; by \cite[Theorem 39]{St}, $v^+ $ is a highest weight vector of $V'$.
Clearly,  $v^+$ does not lie in any maximal $\mathbf{G}$-submodule $M$ of
 $V'$, so $V'/M$ is an \ir $\mathbf{G}$-module of \hw $\om'$.)  By Corollary \ref{no1}(3),
  $\om'\neq \om'_i$ for any $i=1\ld n$.   By Lemma \ref{sd4}(2), either $0\prec \om'$,   or $\om'_1\preceq\om' $ and hence
 $\om'_1\preceq  \mu\preceq \om'$ for every dominant weight $\mu$ of $U$; moreover, as $\om'\neq \om'_i$ for $i=1\ld n$, we have $\om'_1\prec \om'_1+\om'_2\preceq \om'$ (Lemma \ref{sd4}(5)). By Lemma \ref{w0}, $V|_{\mathbf{G}}$ has all weights of the Weyl module $\tilde V_{\om'}$. In turn, all weights of $V_{\om_1'+\om_2'}$
are weights of $\tilde V_{\om'}$, as stated. $\Box$

\subsection{Real elements in $SL_n(q)$ and $SU_n(q)$}

 \bl{zz1}
 Let $ H=SL_n(q)$ or $SU_n(q)$, where $q$ is an arbitrary prime power, and let $g\in H$ be a semisimple element.  Suppose that g is real. Then g is conjugate to a subgroup isomorphic to $Sp_n(q)$ or $Sp_{n-1}(q)$ (depending on parity of n).\el

Proof.   The \mult of every \ei $e$ of $g$ as an element of $GL_n(\overline{\FF}_q)$
equals that of   $e\up$. By determinant reason,  the \mult of $-1$ is even.
If $n$ is odd then 1 is an \ei of $g$ on the natural module $V$ for $H$,  hence $g$ is contained in a subgroup isomorphic to  $SL_{n-1}(q)$ or $SU_{n-1}(q)$.   (In
the unitary case observe that the 1-eigenspace $V^g$ of $g$ on  $V$ is non-degenerate.
So $V^g$ contains an anisotropic vector, and the claim follows again.)  So  $n$ can be assumed to be even. Then the
 result is due to Wall \cite[p.36]{Wa}. (This is straightforwardly for $H=SL_n(q);$
 for $SU_n(q)$ this is explained in \cite[p. 594]{GV}.
 For $n,q$ even the result is explicitly stated
in \cite[Proposition 2.4]{GV}.)  $\Box$

 \bl{yz2}
 Let $h\in H=GL_n(q)$. If $|h|$ divides $q^i+1$ for some integer $i>0$ then $h$ is real.\el

Proof. It suffices to prove this for the case where $h$ is  irreducible and $|h|>2$. Let $E$ be the enveloping
algebra of $\lan h\ran$ over $\FF_q$. By Schur's lemma, $E$ is a field, and  $E\cong \FF_{q^n}$ as $h$ is \irt Then $|h|$ divides $q^n-1$, and hence $(q^n-1,q^i+1)$, and also $(q^n-1,q^{2i}-1)$. By \cite[Hilfsatz 2]{Hu}, $(q^n-1,q^{2i}-1)=q^{(n,2i)}-1$. If $n$ is odd then $q^{(n,2i)}-1=q^{(n,i)}-1=(q^n-1,q^i-1)$, so $|h|$ divides $q^{i}-1$. 
As  $(q^{i}-1,q^{i}+1)\leq 2$, we have $|h|\leq 2$, a contradiction.

Let  $n=2m$ be even and $k=(m,i)$.   Then $(q^n-1,q^i+1)=(q^{2m}-1,q^i+1)$ which divides $(q^{2m}-1,q^{2i}-1)
  =q^{2k}-1$. So $|h|$ divides $q^{2k}-1$. As $2k|n$,
  it follows that $\FF_{q^{2k}}$ is a subfield of $E\cong \FF_{q^n}$, and hence $h$ lies in the subfield
  $E_1$ of $E$ of order $q^{2k}$.   This implies   $E_1=E$ and $2k=n$, $k=m$.
 Then $m|i$ (as $k=(m,i)$), and $i/m$ is odd. Indeed, if $i=2mj=nj$ then $q^i+1=q^{nj}+1$; as $|h|$
 divides $q^n-1$ and hence $q^{nj}-1$, it follows that $|h|\leq 2$, a contradiction. So $i/m$ is odd,
 and hence $q^m+1$ divides $q^i+1$. In addition,    $(q^i+1)/(q^m+1)$ is odd, so $|h|$ divides $q^m+1$.
 (Indeed, $(|h|, q^i-1)\leq 2$, and hence $(|h|, q^m-1)\leq 2$ as $m|i$ and $q^m-1$ divides $q^i-1$. So $|h|$
 divides $2(q^m+1)$. As $|h|$
 divides $q^i+1$ and $q^i+1=a(q^m+1)$ with $a$ odd, we conclude that $|h|$ divides $q^m+1$, as claimed.)

Let $N=N_{H}(E)$. Then the Galois group $\FF_{q^n}/\FF_q$ is isomorphic to a subgroup of $N/C_H(E)$
by the No\"ether-Skolem theorem \cite[\S 12.6]{P}. The Galois group
 $\FF_{q^n}/\FF_{q^m}$ is of order 2, and the non-trivial Galois \au of $\FF_{q^n}/\FF_{q^m}$
 sends $x\in \FF_{q^n}$ to $x^{q^m}$. As $|h|$ divides $q^m+1$, we have $h^{q^m}=h^{q^m+1}h\up=h\up$.
 So $h$ is real in $GL_n(q)$. $\Box$

\bl{ru1}
Let $g\in G=SU_n(q)$. Suppose that  $|g|$ divides $q^i-1$ for some $i$ odd. Then $g$ is real. \el

Proof. It suffices to prove the lemma for the case where $g$ stabilizes no non-degenerate proper subspace of the underlying space $V$ of $G$. Then either $g$ is \ir or $n$ is even and $V=V_1+V_2$, where $V_1,V_2$ are $g$-stable totally isotropic subspaces of dimension $n/2$. In the former case $n$ is odd and  $|g|$ divides $q^n+1$, in the latter case  $g$ is \ir on $V_1$ and  $|g|$ divides $q^n-1$. Note that $(q^i-1,q^j-1)=q^{(i,j)}-1$ by \cite[Hilfsatz 2]{Hu}.

In the former case    $|g|$ divides   $(q^i-1,q^n+1)$, hence also $(q^i-1,q^{2n}-1)=q^{(i,2n)}-1=q^{(i,n)}-1$.   So $|g|$ divides $q^n-1$. As $(q^n-1,q^n+1)\leq 2$, it follows $|g|\leq 2$, and then $g$ is real.

In the latter case  $g\in Y$, where $Y=\{s\in G: sV_1=V_1\}$. It is well known that  $Y\cong GL(V_1)\cong GL_{n/2}(q^2)$. Furthermore, there is a basis $B$, say, in $V$ such that $B\cap V_i$ is a basis of $V_i$ for $i=1,2$, and the matrix of every  $y\in Y$ is $\diag(x,({}^tx^{-1})^J)$ with $x\in GL(V_1)$, where ${}^tx$
means the transpose of $x$ and $J$ is the Galois \au of $\FF_{q^2}/\FF_{q}$ of order 2 extended to $GL_{n/2}(q^2)$.

Let $g=\diag(h,({}^th^{-1})^J))$, so $h$
 is \ir in $GL(V_1)$. 
 If  $h$ is real in $GL(V_1)$ then $g$ is real in $G$, so we assume that $h$ is not real in $GL(V_1)$. In addition, the characteristic \po of $h$ is \ir as well as that of $({}^th^{-1})^J; $ and they are distinct (otherwise $V$ is a homogeneous reducible $F_{q^2}\lan g \ran$-module, and then  $V$ contains a non-degenerate $g$-stable subspace
 (see \cite[Lemma 6.3]{EZ}). Therefore, $C_G(g)V_1=V_1$, so $C_G(g)\subset Y$; as the group $C_{GL(V_1)}(h)$ is cyclic, so is $C_{G}(g)$.

 As $i$ is odd, $n$ is even, we have $(i,n)=(i,n/2)$. Therefore,  $|h|$ divides $(q^i-1,q^n-1)=q^{(i,n)}-1=q^{(i,n/2)}-1$ by the above. As $h$ is \ir on $V_1$, the $\FF_{q^2}$-enveloping algebra $E$, say, of $\lan h\ran$ is a field of order $q^n$, so  the multiplicative group $E^\times$ of $E$ is of order  $q^{n}-1$. The elements of $E^\times$ whose order divides $q^{n/2}-1$ together with 0 form a subfield $E_1$, say, of order $q^{n/2}$. So $h\in E_1$. In addition, $n/2$ is odd (otherwise $E_1$  contains   $\FF_{q^2}$,  which is false as $h
 $ is irreducible). Let $d_1$ be a generator of $E_1^\times$ and
 $d=\diag(d_1,({}^td_1^{-1})^J))$ and $D=\lan d\ran$. Then $g\in D$. 

 Note that $Sp_n(q)$ contains a cyclic subgroup of order $q^{n/2}-1$.
  So the result follows if we show that every cyclic subgroup of $G$ of order $q^{n/2}-1$ is conjugate to $D$.
  For this it suffices to show that $D$ contains a Sylow $r$-subgroup $R$ of $G$ for some prime $r$ such that $C_G(R)=E^\times$. We can assume $n>2$, as $SU_2(q)\cong Sp_2(q)$,
  and $g\in Sp_2(q)$ is real whenever $(|g|,q)=1$. 

 By \cite[Theorem 5.2.14]{KL}, $q^n-1$ is divisible by a prime  $r$  
 such that $r$ does not divide $q^j-1$
 for any $j<n$, unless $q =2,n =6$ or $n= 2$. By the above,  $n> 2$; if $(q,n)=(2,6)$ then we take $r=7$. Then $|R|=7$
  and one can easily check that $C_G(R)=E^\times$. 

 In the non-exceptional case  let $R$ be the Sylow $r$-subgroups of $D$. Then $R$ is
  a Sylow $r$-subgroups of $GL_{n/2}(q^2)$. Comparing the orders of $GL_{n/2}(q^2)$ and $G$ we conclude that
 $R$ remains a Sylow $r$-subgroups of $G$. (Note that $(r,q^j+1)=1$ for  $0<j<n$ as otherwise $r$ divides $q^{2j}-1$
 and $(q^{2j}-1, q^i-1)=q^{(i,j)}-1$ as $i$ is odd; this contradicts the choice of $r$ as $(i,j)<n$.)
 Finally, one easily observe that $C_G(R)=E^\times$. $\Box$

\bl{hg5}
 Let $\mathbf{L}=SL_{n+1}(F)$,  $n$  even, and let $\mathbf{H}=SL_{n+1}(F) $.  Let  $\phi$ be an \irr of $\mathbf{L}$ with $p$-restricted \hw $\lam$. Let    $\lam_1\ld \lam_{n}$ be the fundamental weights of $\mathbf{L}$ and  $\om_1\ld \om_{n-1}$ be those of  $\mathbf{H}$.

$(1)$ \st $n>2$ if $p=2$ and if $\nu$ is the \hw of a composition factor of $\phi|_{\mathbf{H}}$ and $\nu$  is $p$-restricted then $\nu\in\{0,\om_1\ld \om_{n-1}\}$. Then $\lam\in\{\lam_1\ld \lam_n\}$. 

$(2)$ Suppose that $\lam=\lam_i$ with $i\leq n$. Then $\phi_\lam|_{\mathbf{H}}$ contains an \ir constituent of \hw $\om_{i-1}$. (If $i=1$ then $\om_{i-1}$ is interpreted as the zero weight.)

$(3)$ Let $\psi$ be any \rep of  $\mathbf{L}$ and ${\mathbf G}=Sp_{n}(F)$. Then $\psi|_{{\mathbf G}}$ has weight $0$. 
\el

 Proof. We can view $\mathbf{H}$ as a subgroup of $\mathbf{L}$ of shape $\diag(\mathbf{H},1)$ or $\diag(1,\mathbf{H})$.  Then the reference torus $\mathbf{T}_{\mathbf{H}}$ of $\mathbf{H}$ (chosen to be the group of diagonal matrices) is contained in that of $\mathbf{L}$  and a Borel subgroup $B$    of  $\mathbf{H}$ containing $\mathbf{T}_{\mathbf{H}}$ is contained in a Borel subgroup of $\mathbf{L}$ that contains the reference torus $\mathbf{T}$ of $\mathbf{L}$. Then $\ep_n(\mathbf{T}_{\mathbf{H}})=1$ in the former case and $\ep_1(\mathbf{T}_{\mathbf{H}})=1$ in the latter case.

 (1) Let $\lam=\sum a_i\lam_i$. Then the restriction $\phi|_{\mathbf{H}}$ has composition factors with highest weights
 $a_1\om_1+\cdots+a_{n-1}\om_{n-1}$ and $a_2\om_1+\cdots+a_{n}\om_{n-1}$. These weights are $p$-restricted so must lie in $\{0,\om_1\ld \om_{n-1}\}$. This implies the claim unless possibly $\lam=\lam_1+\lam_{n}$.
Then $\phi_{\lam_1+\lam_{n}}$ is well known to be the unique non-trivial \ir constituent of the adjoint module (when $\mathbf{L}$ acts by conjugation of the $((n+1)\times (n+1))$-matrices with trace 0). For $n>3$ one can easily observe that $\phi_{\om_1+\om_{n-1}}$ is a constituent of the adjoint module and hence of  $V_{\lam_1+\lam_{n}}$, whence the result. If $n=2$ then the restriction of
$\phi_{ \lam_1+\lam_2 }$ to $\mathbf{H}$ has an \ir constituent of \hw $2\om_1$, which is $p$-restricted for $p\neq2$. So this
option is rules out.

(2)  $V_\lam$ has weight $\mu:=\ep_1+\cdots +\ep_{i-1}+\ep_{n}$. Let $0\neq v\in V_\lam$ be a vector in the $\mu$-weight space. Then $v$ is a primitive vector for $\mathbf{H}$, that is, $Bv=\lan v\ran$.
Therefore, $V':=\lan \mathbf{H}v\ran$ is a $\mathbf{H}$-module (possibly reducible) with \hw    $\ep_1+\cdots+\ep_{i-1}$ \cite[Theorem 39]{St}. Then one can factorize   $V'$ by a maximal submodule to obtain a composition factor with  highest weight $\om_{i-1}$, as claimed.

(3) It suffices to prove this for tensor-indecomposable \ir representations, and then we can assume that $\psi$ has a $p$-restricted \hw $\mu$, say. Then either $0\prec\mu$ or $\lam_i\preceq\mu$ for some $i\in\{1\ld n\}$. If $i$ is even then $\psi_{\lam_i}|_{\mathbf{G}}$ has weight 0 by Lemma \ref{ww0}. If $i$ is odd then $\psi_{\lam_i}|_{\mathbf {H}}$ has a constituent $\rho$, say, of \hw $\om_{i-1}$; as $i-1$ is even, $\rho|_{\mathbf{G}}$ has weight 0,  by Lemma \ref{ww0}. $\Box$


\begin{corol}\label{rq2} Let 
$H=SL_{n+1}(q)$ or $SU_{n+1}(q)$, where n is even, p a prime dividing q, and let $\rho$ be a p-modular \rep of H.
Then $\rho(g)$ has \ei $1$ for every real element $g\in H$.
\end{corol}

Proof. It suffices to prove the theorem for $\rho$ irreducible, and
with $p$-restricted highest weight.  By Lemma \ref{zz1}, we can assume that $g\in G\cong Sp_{n}(q)$. So the result follows from  Lemma \ref{hg5}(3). $\Box$

 \begin{theo}\label{th5}
  Let $\mathbf{H}=SL_{n+1}(\overline{\FF}_2)$, let  $H=SL_{n+1}(2)$ or $SU_{n+1}(2)$, and let $g\in H$ be a semisimple element. Let $\phi$ be an \irr of $\mathbf{H}$ with highest weight $\lam$.

 $(1)$ Suppose that $g\in H$ is  real.  Then $\phi(g)$ has \ei $1$
unless (possibly)  $n+1$ is even and $\lam=2^i\lam'$, where $\lam' $  is some fundamental weight of $\mathbf{H}$.

   $(2)$ Suppose that $H=SL_{n+1}(2)$ and $|g|$ divides $2^i+1$ for some i,    or $H=SU_{n+1}(2)$  and $|g|$ divides $2^i-1$ for some i odd. Then the conclusion of $(1)$ holds.
 \end{theo}

   Proof. It suffices to prove (1) as $g$ in (2) is real by  Lemmas \ref{yz2} and \ref{ru1}. Furthermore, if $n+1$ is odd then
   the argument in the proof of Corollary \ref{rq2} works and yields the result.
   So we assume that $n+1$ is even. By Lemma \ref{zz1}, we can assume that $g\in G=Sp_{n+1}(2)$.

    Let $\lam_1\ld \lam_n$ be the  fundamental weights of $\mathbf{H}$ and $\om_1\ld\om_{(n+1)/2}$  those   of $\mathbf{G}=Sp_{n+1}({\overline{\FF}_2})$. Suppose first that $\lam$ is 2-restricted.
        By Lemma \ref{do1}, $\rho|_{\mathbf{G}}$ has either weight 0 (and we are done) or all weights of an
        \irr of $\mathbf{G}$ with highest weight  $\om_1+\om_2$.
In the latter case the result follows  by Lemma \ref{ft2}(3).

Suppose that $\lam$ is not 2-restricted. Then $\lam=\sum 2^j\mu_j$, where $\mu_j$ are 2-restricted dominant weights.
Then $\phi=\otimes_j \rho_{2^j\mu_j}$ and $\phi|_H=\otimes_j \rho_{\mu_j}|_H$. Clearly, we only need to examine the case where
this product has at least two terms. If $\mu_j$ is not radical, denote by
$\mu_j'$ the fundamental weight such that $\mu_j'\preceq \mu_j$, otherwise set $\mu_j'=\mu_j$.
Then $\Omega(\otimes_j \rho_{\mu'_j})\subseteq \Omega(\otimes_j \rho_{\mu_j})$, so it suffices to prove the theorem for the case with $\mu_j=\mu'_j$. Set $\nu_j=\mu'_j|_{\mathbf{T}'}$, where $\mathbf{T}'$ is a maximal torus of $\mathbf{G}$ as in Lemma \ref{ww1}, which tells us that $\nu_j\in\{\om_1\ld \om_n\}$ whenever $\mu_j'$ is not radical. For such
$j$ the weights of the Weyl module $\tilde V_{\nu_j}$ are weights of $\rho_{\mu_j}$ (Lemma \ref{w0}). In addition, either
$\nu_j$ is radical or $j$ is odd and $\tilde V_{\nu_j}$ has weights $\pm\om_1$ (Lemmas \ref{sd4} and \ref{pr4} ). So the tensor product of two such modules
has weight 0. \itf that we can write $\phi|_{\mathbf{G}}=\phi_1\otimes \phi_2$, where $\phi_1|_{\mathbf{G}}$ has weight 0
and $\phi_2=\rho_{\mu_j}$ for some $j$, where $\rho_{\mu_j}$ contains all weights of $\tilde V_{\nu_j}$ for some $\nu_j\in\{\om_1\ld \om_n\}$ with $j$ odd. If $\phi_1$ is trivial then the result follows by Corollary \ref{no1}(3). Otherwise,   $\phi_1|_{\mathbf{G}}$ has weight $\om_2$. Indeed, the weights of $\si\in\Omega(\phi_1)$ have the same residue modulo the radical weights (as so are the weights of every multiple $\rho_{\mu_j}$), and  $\si|_{\mathbf{T}'}=0$ for some $\si$. \itf $\si|_{\mathbf{T}'}$ is radical for every $\si\in\Omega(\phi_1)$. In particular,
this is the case for $\si|_{\mathbf{T}'}$ to be the highest weight of a non-trivial composition factor of $\phi_1|_{\mathbf{G}}$. The latter has weight  $\om_2$ by Lemma \ref{pr4}, and hence the weight $\ep_2+\ep_3$ if $n>2$. So all weights of $V_{\om_1+\om_2}$ are weights of $\phi|_{\mathbf{G}}$. 
Then the result follows by Lemma \ref{ft2}.
  $\Box$

\med
{\it Proof of Theorem} \ref{th6}. The result is a special case of Theorem \ref{th5}.

\med
The following Lemmas \ref{co4} and \ref{co5} show that the result of Theorem \ref{th5} is in a sense  best possible.

 \bl{co4}
 Let $h\in H=SL_{2n}(2)$. Suppose that $|h|=p$, where  $p=2^{n}+1$ is a prime.  Let $\phi$ be an $i$-th exterior power of the natural $\FF_2H$-module, where i is odd. Then $\phi(h)$ does not have \ei $1$. \el

 Proof. Note that   $h$ is conjugate to a Singer cycle in $G=Sp_{n}(2)$.   By Lemma \ref{ww0}, the composition factors of $\phi|_{\mathbf{G}}$ are of highest weights $\om_j'$ for $j$ odd. So  Lemma \ref{op1} yields the result. $\Box$

 \bl{co5}
 Let $h\in H=SU_n(2)$. Suppose that $|h|=p$, where  $p=2^{n/2}-1$ is a prime.  Let $\phi$ be an $i$-th exterior power of the natural $\FF_4H$-module, where $i<n-1$ is odd. Then $\phi(h)$ does not have \ei $1$.
 \el

 Proof. This follows from Lemma \ref{016}. Indeed,  $n/2$ is odd, so   $h$ is conjugate to an element of $G=Sp_{n}(2)$ (Lemmas \ref{ru1} and \ref{zz1}).
 Moreover, $h\in G$ is a Singer cycle, so $h$ is a generator of  a (cyclic) maximal torus $T_w$ of $G$ labeled by an element $w$ of the Weyl group which transitively permutes $\ep_1\ld \ep_n$. Let $W_i$ be as in Lemma \ref{016}. Then $w$ is conjugate to an element of $W_i $ \ii $i=n$. As in Lemma \ref{co4} one observes that the composition factors of $\phi|_G$ are of highest weights $\om_j'$ for $j\in\{1\ld n\}$ odd. Then the reasoning in the proof of Lemma \ref{op1} shows
  that $h$ does not have \ei 1 on $V_{\om'_i}$ whenever $i<n-1$ odd. $\Box$

\bigskip
Alexandre E. Zalesski

Department of Physics, Mathematics and Informatics,

National Academy of Sciences of Belarus,

66 Prospekt Nezalejnasti,

220072 Minsk,

Belarus

e-mail: alexandre.zalesski@gmail.com

\end{document}